\documentclass[11pt,twoside]{amsart}
\usepackage{amssymb}

\newcommand{\Z}{{\mathbb{Z}}}
\newcommand{\N}{{\mathbb{N}}}
\newcommand{\R}{{\mathbb{R}}}
\newcommand{\bT}{{\mathbb{T}}}
\newcommand{\bH}{{\mathcal{H}}}
\newcommand{\cA}{{\mathcal{A}}}
\newcommand{\cC}{{\mathcal{C}}}
\newcommand{\cD}{{\mathcal{D}}}
\newcommand{\cE}{{\mathcal{E}}}
\newcommand{\cL}{{\mathcal{L}}}
\newcommand{\cR}{{\mathcal{R}}}
\newcommand{\cLR}{{\mathcal{LR}}}
\newcommand{\cO}{{\mathcal{O}}}
\newcommand{\cS}{{\mathcal{S}}}

\newcommand{\fC}{{\mathfrak{C}}}
\newcommand{\fM}{{\mathcal{M}}}
\newcommand{\fN}{{\mathcal{N}}}
\newcommand{\fR}{{\mathfrak{R}}}
\newcommand{\fS}{{\mathfrak{S}}}
\newcommand{\fI}{{\mathfrak{I}}}
\newcommand{\fX}{{\mathfrak{X}}}
\newcommand{\fs}{{\mathfrak{s}}}
\newcommand{\ft}{{\mathfrak{t}}}
\newcommand{\fp}{{\mathfrak{p}}}
\newcommand{\fu}{{\mathfrak{u}}}
\newcommand{\fv}{{\mathfrak{v}}}
\newcommand{\ba}{{\mathbf{a}}}
\newcommand{\bw}{{\mathbf{w}}}

\newcommand{\ty}{{\tilde{y}}}
\newcommand{\Irr}{{\operatorname{Irr}}}
\newcommand{\Ind}{{\operatorname{Ind}}}

\renewcommand{\leq}{\leqslant}
\renewcommand{\geq}{\geqslant}

\renewcommand{\atop}[2]{\genfrac{}{}{0pt}{}{#1}{#2}}

\newtheorem{thm}{Theorem}[section]
\newtheorem{lem}[thm]{Lemma}
\newtheorem{cor}[thm]{Corollary}
\newtheorem{prop}[thm]{Proposition}

\theoremstyle{definition}
\newtheorem{defn}[thm]{Definition}
\newtheorem{exmp}[thm]{Example}

\theoremstyle{remark}
\newtheorem{rem}[thm]{Remark}

\begin{document}

\title{Kazhdan--Lusztig cells and the Murphy basis}

\author{Meinolf Geck}
\address{Institut Girard Desargues, bat. Jean Braconnier, Universit\'e Lyon 1,
21 av Claude Bernard, F--69622 Villeurbanne Cedex, France}
\email{geck@igd.univ-lyon1.fr}
\date{April 6, 2005}

\begin{abstract} Let $\bH$ be the Iwahori--Hecke algebra associated with
$\fS_n$, the symmetric group on $n$ symbols. This algebra has two important 
bases: the Kazhdan--Lusztig basis and the Murphy basis. While the former
admits a deep geometric interpretation, the latter leads to a purely 
combinatorial construction of the representations of $\bH$, including 
the Dipper--James theory of Specht modules. In this paper, we establish a 
precise connection between the two bases, allowing us to give, for the 
first time, purely algebraic proofs for a number of fundamental properties 
of the Kazhdan--Lusztig basis and Lusztig's results on the $\ba$-function. 
\end{abstract}

\maketitle

\pagestyle{myheadings}

\markboth{Geck}{Kazhdan--Lusztig cells and the Murphy basis}

\section{Introduction} \label{sec0}
Let $W$ be a finite or affine Weyl group and $\bH$ be the associated generic 
Iwahori--Hecke algebra. By definition, $\bH$ is equipped with a standard 
basis usually denoted by $\{T_w\mid w\in W\}$. In a fundamental paper, 
Kazhdan and Lusztig \cite{KaLu} constructed a new basis $\{C_w\mid w\in W\}$
and used this, among other applications, to define representations of $\bH$ 
endowed with canonical bases. While that construction is completely 
elementary, it has deep connections with the geometry of flag varieties 
and the representation theory of Lie algebras and groups of Lie type. An 
excellent account of these connections is given in Lusztig's survey
\cite{Lusztig90b}. 

One of the main consequences of the geometric interpretation are certain 
``positivity properties'' for which no elementary proofs have ever been 
found. In turn, these positivity properties allowed Lusztig to establish a 
number of fundamental properties of the basis $\{C_w\}$, which are concisely 
summarized in a list of $15$ items (P1)--(P15) in \cite[14.2]{Lusztig03}. 
For example, one application of these properties is the construction of
a ``canonical'' isomorphism from $\bH$ onto the group algebra of $W$ (when
$W$ is finite); see Lusztig \cite{Lu0}. 

This paper arose from the desire to give elementary, purely algebraic proofs 
for (P1)--(P15) in the case where $W$ is the symmetric group $\fS_n$. This
goal will be achieved in Secetion~5; we will also prove a tiny piece of 
``positivity'' in this case, namely, the fact that the structure constants 
of Lusztig's ring $J$ are $0$ or $1$. The key idea of our approach is to 
relate the basis $\{C_w\}$ to the basis constructed by Murphy \cite{Mu1}, 
\cite{Mu2}. The main problem in establishing that relation is that the 
Murphy basis elements are not directly indexed by the elements of $W=\fS_n$. 
Indeed, those elements are written as $x_{\fs\ft}$ where $(\fs,\ft)$ runs 
over all pairs of standard $\lambda$-tableaux, for various partitions 
$\lambda$ of $n$. Now, the Robinson--Schensted correspondence does associate 
to each element of $\fS_n$ a pair of standard tableaux of the same shape, 
but this works on a purely combinatorial level; we need to relate this to 
basis elements of $\bH$. Eventually, we shall see that the ``leading matrix 
coefficients'' introduced by the author \cite{my02} provide a bridge to 
pass from the Kazhdan--Lusztig basis to the Murphy basis of $\bH$; see 
Theorem~\ref{ourmain} and Corollary~\ref{bridge}. While the explicit form 
of the base change seems to be rather complicated, our results are 
sufficiently fine to enable us to translate combinatorial properties of 
the Murphy basis into properties of the Kazhdan--Lusztig basis.

This paper is organized as follows. In Section~2, we recall the basic 
definitions concerning the Kazhdan--Lusztig basis of $\bH$, the left cells 
and the corresponding representations. In Section~3, we consider the 
orthogonal representations and leading matrix coefficients introduced by the 
author \cite{my02}. These provide an efficient tool for identifying 
Kazhdan--Lusztig basis elements and irreducible representations occuring as 
constituents in left cell representations. 

In the remaining parts of this paper, we exclusively consider the case 
where $W=\fS_n$ is the symmetric group. In Section~4, we recall the 
construction of the Murphy basis $\{x_{\fs\ft}\}$ of $\bH$. One of 
Murphy's key results is that, for a fixed partition $\lambda$, the 
submodule $N^\lambda$ spanned by all basis elements $x_{\fs\ft}$, 
where $(\fs,\ft)$ are standard $\mu$-tableaux such that $\lambda 
\trianglelefteq \mu$, is a two-sided ideal in $\bH$. In 
Theorem~\ref{ourmain}, we show that the image of $N^\lambda$ under a 
certain algebra automorphism is spanned by the Kazhdan--Lusztig basis 
elements $C_w$, where $w$ has a non-zero leading matrix coefficient in a 
representation labelled by a partition $\mu$ such that $\lambda 
\trianglelefteq \mu$. 

Finally, in Section~5, we discuss the applications to the left cells 
in $\fS_n$. In Theorem~\ref{ourmain1}, we characterize the Kazhdan--Lusztig
pre-order relation $\leq_{\cLR}$ in terms of the dominance order on
partitions. In Theorem~\ref{ourmain2}, we establish Lusztig's 
``Property (A)'' which played a decisive role in \cite{Lu0}. In 
Theorem~\ref{finalrs}, we interprete our results in terms of the 
Robinson--Schensted correspondence; our approach even yields a new 
proof for a key step in establishing the fact that the left cells of
$\fS_n$ are given by the Robinson--Schensted correspondence. To deal with 
Lusztig's $\ba$-function, we rely on the methods developped by Iancu and 
the author in \cite[\S 4]{GeIa05}. These methods show that most of the 
properties (P1)--(P15) follow from a relatively small set of hypotheses 
which are easily seen to be satisfied thanks to the link between the 
Kazhdan--Lusztig basis and the Murphy basis; see Theorem~\ref{p114}.

\section{Kazhdan-Lusztig cells} \label{sec-kl}

Let $W$ be a Coxeter group with (finite) generating set $S$. Until the
end of this section, we do not need to assume that $W$ is finite. Let  
$l\colon W \rightarrow \N$ be the usual length function with respect to 
$S$ (where $\N=\{0,1,2, \ldots\}$). Let $A={\Z}[v,v^{-1}]$ be the ring of 
Laurent polynomials in an indeterminate~$v$. Let $\bH=\bH_A(W,S)$ be the 
corresponding Iwahori--Hecke algebra. Then $\bH$ is free over $A$ with 
basis $\{T_w\mid w\in W\}$; the multiplication is given by the rule
\[ T_sT_w=\left\{\begin{array}{cl} T_{sw} & \quad \mbox{if $l(sw)=l(w)+1$},
\\ T_{sw}+(v-v^{-1})T_w & \quad \mbox{if $l(sw)=l(w)-1$},\end{array}
\right.\]
where $w\in W$ and $s\in S$.  For basic properties of $W$ and $\bH$, we
refer to \cite{ourbuch}.

\subsection{Kazhdan--Lusztig bases} \label{klbase} There are two types of
Kazhdan--Lusztig bases, denoted by $\{C_w\mid w\in W\}$ and $\{C_w'\mid 
w\in W\}$ in the original article by Kazhdan and Lusztig \cite{KaLu}. 
In \cite{Lusztig03}, Lusztig writes $c_w$ instead of $C_w'$, but we shall 
adhere to the older notation here. The precise definitions are as follows.

There is a unique ring involution $A\rightarrow A$, $a \mapsto \bar{a}$, 
such that $\bar{v}=v^{-1}$. We can extend this map to a ring involution 
$j \colon \bH \rightarrow \bH$ such that 
\[ j\bigl(\sum_{w \in W} a_w T_w\bigr)=\sum_{w \in W} \varepsilon_w\,
\bar{a}_w T_w \qquad (a_w \in A),\]
where we write $\varepsilon_w=(-1)^{l(w)}$ for any $w\in W$.
Furthermore, there is an $A$-algebra automorphism
\[ \dagger\colon \bH\rightarrow \bH,\qquad T_s \mapsto T_s^\dagger=-T_s^{-1} 
\quad (s \in S).\]
It is easily checked that $j$ and $\dagger$ commute with each other. For 
any $h\in \bH$, we set $\overline{h}:=j(h)^\dagger=j(h^\dagger)$. Thus, the 
map $\bH \rightarrow \bH$, $h \mapsto \overline{h}$, is a ring involution 
such that
\[ \overline{\sum_{w \in W} a_w T_w}=\sum_{w \in W} \bar{a}_w
T_{w^{-1}}^{-1} \qquad (a_w \in A).\]
Now, for each $w\in W$, there exists a unique element $C_w'\in \bH$
such that 
\[\overline{C}_w'=C_w'\qquad \mbox{and}\qquad C_w'\equiv T_w \quad\bmod
\bH_{<0},\]
where $\bH_{<0}:=\sum_{w\in W} A_{<0}\, T_w$ and 
$A_{<0}:=v^{-1}{\Z}[v^{-1}]$; see \cite[Theorem~5.2]{Lusztig03}. The 
elements $\{C_w'\mid w\in W\}$ form an $A$-basis of $\bH$, and we have 
\[ C_w'=T_w+\sum_{\atop{y \in W}{y < w}} p_{y,w} \,T_y, \]
where $<$ denotes the Bruhat--Chevalley order on $W$ and $p_{y,w}\in
A_{<0}$ for all $y<w$ in $W$. It will be technically more convenient to 
work with the $C$-basis of $\bH$. (The reasons can be seen, for example, 
in \cite[Chap.~18]{Lusztig03}.) We set $C_w:=\varepsilon_w j(C_w')$. Then 
we have
\[ C_w=T_w+\sum_{\atop{y\in W}{y<w}} \varepsilon_y \,\varepsilon_w \,
\,\overline{p}_{y,w}\,T_y  \qquad \mbox{for all $w\in W$}.\]
As before, one shows that the basis element $C_w$ is uniquely determined
by the conditions that
\[\overline{C}_w=C_w\qquad \mbox{and}\qquad C_w\equiv T_w \quad\bmod
\bH_{>0},\]
where $\bH_{>0}:=\sum_{w\in W} A_{>0}\, T_w$ and $A_{>0}:=v{\Z}[v]$. 

\subsection{Multiplication rules} \label{multrule}
For any $x,y\in W$, we write
\[ C_x'\,C_y'=\sum_{z\in W} h_{x,y,z} \, C_z' \qquad \mbox{where
$h_{x,y,z} \in A$}.\]
We have the following more explicit formula for $s\in S$, $y\in W$
(see \cite[\S 6]{Lusztig03}):
\[ C_s'\,C_y' = \left\{\begin{array}{ll} \displaystyle{C_{sy}'+
\sum_{\atop{z\in W}{sz<z<y}} \mu_{z,y} C_z'} &\quad \mbox{if $sy>y$},\\
(v+v^{-1})\,C_y' &\quad \mbox{if $sy<y$},
\end{array}\right.\]
where $\mu_{z,y}\in \Z$ is the coefficient of $v^{-1}$ in $p_{z,y}$;
see \cite[Cor.~6.5]{Lusztig03}. Using the relation $C_w=\varepsilon_w\,
j(C_w')$, we obtain the formula
\[ C_x\, C_y=\sum_{z\in W} \varepsilon_x\,\varepsilon_y\,\varepsilon_z\,
h_{x,y,z}\, C_z.\]
Note also that, for $s\in S$, we have $C_s'=T_s+v^{-1}T_1$ and $C_s=T_s-
vT_1$. 

\subsection{The Kazhdan--Lusztig pre-orders} \label{klpre}
As in \cite[\S 8]{Lusztig03}, we write $x \leftarrow_{\cL} y$ if 
there exists some $s\in S$ such that $h_{s,y,x}\neq 0$, that is, 
$C_x'$ occurs in $C_s'\, C_y'$ (when expressed in the $C'$-basis) 
or, equivalently, $C_x$ occurs in $C_s\, C_y$ (when expressed in the 
$C$-basis). The Kazhdan--Lusztig left pre-order $\leq_{\cL}$ is the 
relation on $W$ generated by $\leftarrow_{\cL}$, that is, we have 
$x\leq_{\cL} y$ if there exists a sequence $x=x_0,x_1, \ldots,x_k=y$ 
of elements in $W$ such that $x_{i-1} \leftarrow_{\cL} x_i$ for all~$i$. 
The equivalence relation associated with $\leq_{\cL}$ will be denoted 
by $\sim_{\cL}$ and the corresponding equivalence classes are called 
the {\em left cells} of $W$.
                                                                            
Similarly, we can define a pre-order $\leq_{\cR}$ by considering
multiplication by $C_s'$ on the right in the defining relation. The
equivalence relation associated with $\leq_{\cR}$ will be denoted by
$\sim_{\cR}$ and the corresponding equivalence classes are called the
{\em right cells} of $W$.  We have
\[ x \leq_{\cR} y \quad \Leftrightarrow \quad x^{-1} \leq_{\cL} y^{-1}.\]
This follows by using the antiautomorphism $\flat\colon \bH\rightarrow
\bH$ given by $T_w^\flat=T_{w^{-1}}$; we have $C_w'^\flat=C_{w^{-1}}'$
and $C_w^\flat=C_{w^{-1}}$ for all $w\in W$; see \cite[5.6]{Lusztig03}. 
Thus, any statement concerning the left pre-order relation $\leq_{\cL}$ 
has an equivalent version for the right pre-order relation $\leq_{\cR}$, 
via $\flat$. 

Finally, we define a pre-order $\leq_{\cLR}$ by the condition that 
$x\leq_{\cLR} y$ if there exists a sequence $x=x_0,x_1, \ldots, x_k=y$ 
such that, for each $i \in \{1,\ldots,k\}$, we have $x_{i-1} \leq_{\cL} 
x_i$ or $x_{i-1}\leq_{\cR} x_i$. The equivalence relation associated with 
$\leq_{\cLR}$ will be denoted by $\sim_{\cLR}$ and the corresponding 
equivalence classes are called the {\em two-sided cells} of $W$.
                                                                            
\subsection{Lusztig's conjectures} \label{conj}
For the convenience of the reader, we recall here Lusztig's conjectures
(P1)--(P15) from \cite[Chap.~14]{Lusztig03}. (In fact, Lusztig formulates
these conjectures in a slightly more general setting, where $\bH$ is 
defined with respect to a weight function $L$ on $W$; but we will not need 
to discuss that more general situation here.) For a fixed $z\in W$, we set
\[ \ba(z):= \min \{i\geq 0\mid v^i\,h_{x,y,z} \in {\Z}[v]\mbox{ for 
all $x,y\in W$}\};\]
this is Lusztig's function $\ba \colon W \rightarrow \N$, introduced in 
\cite{Lu1}. Note that, for infinite $W$, it is not at all clear if  
$\ba(z)<\infty$ for all $z\in W$. However, Lusztig \cite[13.4]{Lusztig03} 
conjectures that this is always the case. 

Furthermore, we define  
$\Delta(z)\in \Z$ and $0\neq n_z\in \Z$ by the condition
\[p_{1,z}=n_z\,v^{-\Delta(z)}+\mbox{combination of smaller powers of
$v$};\]
note that $\Delta(z)\geq 0$. Given $x,y,z\in W$, we define 
$\gamma_{x,y,z^{-1}}\in \Z$ by
\[ \gamma_{x,y,z^{-1}}=\mbox{constant term of $v^{\ba(z)}\, h_{x,y,z} 
\in {\Z}[v]$}.\]
Let 
\[\cD:=\{z\in W \mid \ba(z)=\Delta(z)\}.\]  
Then Lusztig \cite[14.2]{Lusztig03} conjectures that the following properties
hold. 
\begin{itemize}
\item[\bf P1.] For any $z\in W$ we have $\ba(z)\leq \Delta(z)$.
\item[\bf P2.] If $d \in \cD$ and $x,y\in W$ satisfy $\gamma_{x,y,d}\neq 0$,
then $x=y^{-1}$.
\item[\bf P3.] If $y\in W$, there exists a unique $d\in \cD$ such that
$\gamma_{y^{-1},y,d}\neq 0$.
\item[\bf P4.] If $z'\leq_{\cL\cR} z$ then $\ba(z')\geq \ba(z)$. Hence, if
$z'\sim_{\cL\cR} z$, then $\ba(z)=\ba(z')$.
\item[\bf P5.] If $d\in \cD$, $y\in W$, $\gamma_{y^{-1},y,d}\neq 0$, then
$\gamma_{y^{-1},y,d}=n_d=\pm 1$.
\item[\bf P6.] If $d\in \cD$, then $d^2=1$.
\item[\bf P7.] For any $x,y,z\in W$, we have $\gamma_{x,y,z}=\gamma_{y,z,x}$.
\item[\bf P8.] Let $x,y,z\in W$ be such that $\gamma_{x,y,z}\neq 0$. Then
$x\sim_{\cL} y^{-1}$, $y \sim_{\cL} z^{-1}$, $z\sim_{\cL} x^{-1}$.
\item[\bf P9.] If $z'\leq_{\cL} z$ and $\ba(z')=\ba(z)$, then $z'\sim_{\cL}z$.
\item[\bf P10.] If $z'\leq_{\cR} z$ and $\ba(z')=\ba(z)$, then $z'\sim_{\cR}z$.
\item[\bf P11.] If $z'\leq_{\cL\cR} z$ and $\ba(z')=\ba(z)$, then
$z'\sim_{\cL\cR}z$.
\item[\bf P12.] Let $I\subset S$ and $W_I$ be the parabolic subgroup
generated by $I$. If $y\in W_I$, then $\ba(y)$ computed in terms of $W_I$
is equal to $\ba(y)$ computed in terms of~$W$.
\item[\bf P13.] Any left cell $\fC$ of $W$ contains a unique element
$d\in \cD$. We have $\gamma_{x^{-1},x,d}\neq 0$ for all $x\in \fC$.
\item[\bf P14.] For any $z\in W$, we have $z \sim_{\cL\cR} z^{-1}$.
\item[\bf P15.] Let $\breve{v}$ be a second indeterminate and let 
$\breve{h}_{x,y,z}\in {\Z}[\breve{v},\breve{v}^{-1}]$ be obtained from
$h_{x,y,z}$ by the substitution $v\mapsto \breve{v}$. If $x,x',y,w\in W$ 
satisfy $\ba(w)=\ba(y)$, then
\[\sum_{y' \in W} \breve{h}_{w,x',y'}\,h_{x,y',y}=
\sum_{y'\in W} h_{x,w,y'}\, \breve{h}_{y',x',y}.\]
\end{itemize}
If $W$ is a finite or affine Weyl group then all these properties hold 
thanks to a geometric interpretation which yields the following  
``positivity property'':
\[ p_{x,y}\in v^{-1}{\N}[v^{-1}] \qquad \mbox{and}\qquad h_{x,y,z}\in 
{\N}[v,v^{-1}]\qquad (x,y,z\in W);\]
see the discussion of the ``split case'' by Lusztig 
\cite[Chap.~15]{Lusztig03} and the references there; see also 
Springer \cite{Spr}. 

Once (P1)--(P15) are known to hold, one can construct the ring $J$ as 
explained in \cite[Chap.~18]{Lusztig03}. As an abelian group, $J$ is free 
with a basis $\{t_w\mid w\in W\}$. The multiplication is given by
\[ t_x \cdot t_y=\sum_{z\in W} \gamma_{x,y,z^{-1}} \, t_z\qquad
\mbox{for all $x,y \in W$},\]
where the identity element is $1_J=\sum_{d\in \cD} n_dt_d$. Furthermore, 
we have an $A$-algebra homomorphism $\phi \colon \bH\rightarrow A 
\otimes_{\Z} J$ defined by Lusztig's formula 
\cite[Theorem~18.9]{Lusztig03}. Then all the methods developped in 
\cite[Chap.~20--24]{Lusztig03} can be applied to the study of the
left cell representations of $\bH$.

\subsection{Left cell representations} \label{leftrep} Let $\fC$ be a 
left cell or, more generally, a union of left cells of $W$. We define an 
$\bH$-module by $[\fC]_A:={\fI}_{\fC}/\hat{\fI}_{\fC}$, where
\begin{align*}
{\fI}_{\fC} &:=\langle C_w\mid w\leq_{\cL}
z\mbox{ for some $z \in\fC$}\rangle_A,\\
\hat{\fI}_{\fC} &:=\langle C_w\mid w \not\in \fC, w\leq_{\cL} z
\mbox{ for some $z \in\fC$}\rangle_A.
\end{align*}
Note that, by the definition of the pre-order relation $\leq_{\cL}$, these
are left ideals in $\bH$. Denote by $b_x$ ($x\in \fC$) the residue
class of $\varepsilon_x C_x$ in $[\fC]_A$. Then the elements $\{b_x\mid x
\in \fC\}$ form an $A$-basis of $[\fC]_A$ and the action of $C_w$ ($w\in W$) 
is given by the formula
\[ C_w.b_x=\varepsilon_w\sum_{y \in \fC} h_{w,x,y}\, b_y.\]
Note that, since $C_w^\dagger=\varepsilon_w\, C_w'$, this coincides with 
the definition in \cite[\S 21.1]{Lusztig03}. Up to the change of basis
$b_x \leftrightarrow \varepsilon_x b_x$ ($x \in \fC$), this also 
coincides with the original definition by Kazhdan--Lusztig \cite{KaLu}.

\subsection{Kazhdan--Lusztig's star operation} \label{star} Let us consider
two generators $s \neq t$ in $S$ such that $sts=tst$. We set 
\[ \cD_R(s,t):=\{w\in W\mid \mbox{either $ws<w$, $wt>w$ or 
$wt<w$, $ws>w$}\}.\]
If $w\in \cD_R(s,t)$, then exactly one of the elements $ws,wt$ belongs to
$\cD_R(s,t)$; we denote it $w^*$. The map
\[ \cD_R(s,t) \rightarrow \cD_R(s,t), \qquad  w\mapsto  w^*,\]
is an involution. It is readily checked that  
\begin{equation*}
w \sim_{\cR} w^* \qquad \mbox{for any $w\in \cD_R(s,t)$}.\tag{a}
\end{equation*}
Now let $w,w_1\in W$. Following Kazhdan--Lusztig \cite[\S 5]{KaLu}, we 
write $w \approx w_1$ if there exist $s,t\in S$ as above such that 
$w\in \cD_R(s,t)$ and $w_1=w^*$. Now, the relations $\approx$ and 
$\sim_{\cL}$ are compatible, in the following sense. Let $y\in W$ be such 
that $y\sim_{\cL} w$. Then, by \cite[8.6]{Lusztig03}, we also have $y\in 
\cD_R(s,t)$. Hence, by \cite[Cor.~4.3]{KaLu}, we have $y^* \sim_{\cL} w^*$. 
Consequently, we have a bijection 
\begin{equation*}
\fC \stackrel{\sim}\rightarrow \fC_1, \qquad x\mapsto x_1:=x^*,\tag{b}
\end{equation*}
where $\fC$ is the left cell containing $w$ and $\fC_1$ is the left
cell containing $w_1$. We shall also write $\fC \approx \fC_1$ in this
situation. By \cite[Theorem~4.2]{KaLu} (see also the discussion in 
\cite[\S 5]{KaLu}), the above bijection has the following property:
\begin{equation*}
h_{s,x,y}=h_{s,x_1,y_1} \qquad \mbox{for all $s\in S$ and all $x,y\in \fC$}.
\tag{c}
\end{equation*}
This means that the action of $C_s$ ($s \in S$) on the $\bH$-modules 
$[\fC]_A$ and $[\fC_1]_A$ is given by exactly the same formulas with 
respect to the standard bases of $[\fC]_A$ and $[\fC_1]_A$, respectively. 
Since the elements $C_s$ ($s \in S$) generate $\bH$ as an $A$-algebra,
we can even conclude that 
\begin{equation*}
h_{w,x,y}=h_{w,x_1,y_1} \qquad \mbox{for all $w\in W$ and all $x,y\in \fC$}.
\tag{c'}
\end{equation*}
We also note the following property of the bijection in (b):
\begin{equation*}
l(x)+l(y) \equiv l(x_1)+l(y_1) \quad \bmod 2, \qquad \mbox{for 
all $x, y \in \fC$}.\tag{d}
\end{equation*}
Indeed, by the definition of $w^*$, we have $l(w^*)=l(w)\pm 1$. This
immediately yields (d).

\subsection{Induction from parabolic subgroups} \label{indu}
Let $I\subseteq S$ and consider the parabolic subgroup $W_I\subseteq W$. 
Let $X_I$ be the set of distinguished left coset representatives; we have
\[ X_I=\{w \in W \mid \mbox{$w$ has minimal length in $wW_I$}\}.\]
Furthermore, the map $X_I \times W_I \rightarrow W$, $(x,u) \mapsto xu$,
is a bijection and we have $l(xu)=l(x)+l(u)$ for all $u\in W_I$ and all
$x\in X_I$; see \cite[\S 2.1]{ourbuch}. Let $\bH_I=\langle T_w \mid w \in 
W_I\rangle_A\subseteq \bH$ be the corresponding parabolic subalgebra 
of $\bH$. For any $w\in W_I$, we have $C_w\in \bH_I$ and $C_w'\in \bH_I$; 
hence $\{C_w \mid w\in W_I\}$ and $\{C_w' \mid w\in W_I\}$ are the 
Kazhdan--Lusztig bases of $\bH_I$. 

The linear map $\varepsilon_I \colon \bH_I \rightarrow A$ defined by
$\varepsilon_I(T_w)=\varepsilon_w v^{-l(w)}$ for any $w\in W_I$ is
an algebra homomorphism, called the sign representation. We denote by 
\[ \Ind_I^S(\varepsilon_I),\] 
the $\bH$-module obtained by induction from $\varepsilon_I$. According
to the situation, we will also use the same symbol for the corresponding
character. See \cite[\S 9.1]{ourbuch} for the definition and 
basic properties of induced modules.

\addtocounter{thm}{7}
\begin{lem} \label{cwi} Assume that $W_I$ is finite and let $w_I\in W_I$ 
be the unique element of maximal length. Then the following hold.
\begin{itemize}
\item[(a)] For any $w\in W_I$, we have $T_wC_{w_I}=\varepsilon_w \,v^{-l(w)}
C_{w_I}$. 
\item[(b)] We have $C_{w_I}^2=\varepsilon_{w_I}v^{-l(w_I)} P_I C_{w_I}$, where 
$P_I= \sum_{w\in W_I} v^{2l(w)}$. 
\item[(c)] The set $X_I w_I$ is a union of left cells in $W$; we have
\[ X_I w_I=\{w\in W \mid w\leq_{\cL} w_I\}.\]
We have $[X_I w_I]_A \cong \Ind_I^S(\varepsilon_I)\cong \bH C_{w_I}$ 
(isomorphisms as left $\bH$-modules).
\end{itemize}
\end{lem}

\begin{proof} (a) It is well-known that $l(ww_I)=l(w_I)-l(w)$ for all
$w\in W_I$. Hence, by \cite[Cor.~12.2]{Lusztig03}, we obtain
\[C_{w_I}'=\sum_{w\in W_I} v^{-l(ww_I)}\, T_w =v^{-l(w_I)} 
\sum_{w\in W_I} v^{l(w)}\,T_w.\]
Applying $j\colon \bH \rightarrow \bH$ yields the expression:
\[C_{w_I}=\varepsilon_{w_I}v^{l(w_I)}\sum_{w\in W_I} \varepsilon_w v^{-l(w)}
T_w.\]
Now let $s\in I$. Then we have $sw_I<w_I$ and so the multiplication 
rule in (\ref{multrule}) shows that $C_sC_{w_I}=-(v+v^{-1})C_{w_I}$. Since 
$C_s=T_s-vT_1$, this yields $T_sC_{w_I}=-v^{-1}C_{w_I}$. The required formula 
for $T_wC_{w_I}$ now follows by a simple induction on $l(w)$.

(b) Using (a), we have 
\begin{align*} 
C_{w_I}^2&=\varepsilon_{w_I}v^{l(w_I)}\sum_{w\in W_I} \varepsilon_w 
v^{-l(w)} T_w C_{w_I}\\
&=\varepsilon_{w_I}v^{l(w_I)}\sum_{w\in W_I} v^{-2l(w)} C_{w_I}\\
&=\varepsilon_{w_I}v^{-l(w_I)}\sum_{w\in W_I} v^{2(l(w_I)-l(w))} C_{w_I}\\
&= \varepsilon_{w_I}v^{-l(w_I)}P_I \,C_{w_I}.
\end{align*}

(c) Let $w\in W$ be such that $w\leq_{\cL} w_I$. Then the right descent set
of $w_I$ (which is $I$) is contained in the right descent set of $w$; see
\cite[Lemma~8.6]{Lusztig03}. Hence we can write $w=xw_I$ where $x\in X_I$.
Conversely, if $x\in X_I$ then $l(xw_I)=l(x)+l(w_I)$ and so $xw_I\leq_{\cL} 
w_I$.  This yields the equality $X_Iw_I=\{w\in W\mid w\leq_{\cL} w_I\}$. That 
equality also shows that $X_Iw_I$ is a union of left cells. Hence the module 
$[X_Iw_I]_A$ is defined. Now it is known that $C_{xw_I}$ is the sum of $T_x
C_{w_I}$ and a linear combination of terms $T_y C_{w_I}$ where $y \in X_I$ 
and $y<x$. (This is just a special case of \cite[Prop.~3.3]{myind}, for
example.) Hence we also have that $T_xC_{w_I}$ is the sum of $C_{xw_I}$ and 
a linear combination of terms $C_{yw_I}$ where $y\in X_I$ and $y<x$. 
Consequently, we have 
\[ \langle  C_{xw_I}\mid x\in X_I\rangle_A=\langle T_xC_{w_I} \mid x\in 
X_I\rangle_A,\]
and it is easily seen that the subspace on the right hand side equals 
$\bH C_{w_I}$.  Hence, by the definition of induced modules, we have 
\[ [X_Iw_I]_A \cong \Ind_I^S(\langle C_{w_I} \rangle_A) \cong \bH\,C_{w_I},\]
where $\langle C_{w_I}\rangle_A$ is an $\bH_I$-submodule of $\bH_I$
affording $\varepsilon_I$.
\end{proof}

\section{Orthogonal representations} \label{sec-orth}

We now recall the basic facts concerning the leading matrix coefficients
introduced in \cite{my02}. For this purpose, we assume from now on that
$W$ is a {\em finite} Coxeter group. We extend scalars from $A$ to the 
field $K={\R}(v)$. Every element $x\in K$ can be written in the form
\[ x=r_x\,v^{\gamma_x}f/g\qquad \mbox{where $r_x \in \R$, $\gamma_x \in
\Z$ and $f,g\in 1+v{\R}[v]$};\]
note that, if $x\neq 0$, then $r_x$ and $\gamma_x$ indeed are
{\em uniquely determined} by $x$; for $x=0$, we have $r_0=0$ and we set
$\gamma_0:= +\infty$ by convention. Let 
\[{\cO}:=\{x \in K \mid \gamma_x \geq 0\} \qquad \mbox{and}\qquad
{\fp}:=\{x \in K \mid \gamma_x >0\}.\]
Note that $\cO$ is nothing but the localisation of ${\R}[v]$ in the ideal
$(v)$; hence, $\cO$ is a discrete valuation ring and $\fp$ is the unique 
maximal ideal of $\cO$. The group of units in $\cO$ is given by
\[ {\cO}^\times=\{x \in \cO \mid r_x\neq 0, \gamma_x=0\}.\]
Note that we have
\[ \cO \cap {\R}[v,v^{-1}]={\R}[v] \qquad \mbox{and}\qquad
\fp\cap {\R}[v,v^{-1}]=v{\R}[v].\]
The substitution $v\mapsto 0$ defines an $\R$-linear ring homomorphism 
$\cO \rightarrow \R$ with kernel $\fp$. The image of $x\in \cO$ in $\R$ 
is called the {\em constant term} of $x$. Thus, the constant term of $x$ 
is $0$ if $x\in \fp$; the constant term equals $r_x$ if $x\in \cO^\times$.

\subsection{Schur elements and $a$-invariants} \label{symmtrace}
Extending scalars from $A$ to $K$, we obtain a finite dimensional 
$K$-algebra $\bH_K=K\otimes_A \bH$, with basis $\{T_w\mid w\in W\}$. It is 
well-known that $\bH_K$ is split semisimple and abstractly isomorphic to the 
group algebra of $W$ over $K$; see, for example, \cite[Theorem~9.3.5]{ourbuch}.
Let $\Irr(\bH_K)$ be the set of irreducible characters of $\bH_K$. We write 
this set in the form
\[ \Irr(\bH_K)=\{\chi_\lambda \mid \lambda\in \Lambda\},\]
where $\Lambda$ is some finite indexing set. We have a symmetrizing trace
$\tau \colon \bH_K \rightarrow K$ defined by $\tau(T_1)=1$ and $\tau(T_w)=0$ 
for $1\neq w\in W$. We have
\[ \tau(T_wT_{w'})=\left\{\begin{array}{cl} 1 & \qquad
\mbox{if $w'=w^{-1}$},\\ 0 & \qquad \mbox{if $w' \neq w^{-1}$};
\end{array}\right.\]
see \cite[\S 8.1]{ourbuch}. The fact that $\bH_K$ is split semisimple yields 
that
\[ \tau=\sum_{\lambda\in \Lambda} \frac{1}{c_\lambda} \, \chi_\lambda
\qquad \mbox{where $0 \neq c_\lambda \in {\R}[v,v^{-1}]$}.\]
The elements $c_\lambda$ are called the {\em Schur elements}. By
\cite[8.1.8]{ourbuch}, we have $c_\lambda=P_W/D_\lambda$ where
$P_W=\sum_{w\in W} v^{2l(w)}$ is the Poincar\'e polynomial of $W$ and 
$D_\lambda$ is the ``generic degree'' associated with $\chi_\lambda$. 
There is a unique $a(\lambda)\in \N$ and a positive real number $r_\lambda$
such that 
\[ c_\lambda=r_\lambda \, v^{-2a(\lambda)}+\mbox{ combination of higher 
powers of $v$};\]
see \cite[Chap.~4]{LuBook} and \cite[Def.~3.3]{my02}. The number 
$a(\lambda)$ is called the $a$-invariant of $\chi_\lambda$. 

\subsection{Orthogonal representations} \label{orthrep}
By \cite[Prop.~4.3]{my02}, every $\chi_\lambda$ is afforded by a so-called
orthogonal representation.  This means that there exists a matrix
representation $\fX_\lambda \colon \bH_K\rightarrow M_{d_\lambda}(K)$ with
character $\chi_\lambda$ and an invertible diagonal matrix $P \in
M_{d_\lambda}(K)$ such that the following conditions hold:
\begin{itemize}
\item[(O1)] We have $\fX_\lambda(T_{w^{-1}})=P^{-1} \cdot \fX_\lambda
(T_w)^{\text{tr}}\cdot P$ for all $w\in W$, and
\item[(O2)] the diagonal entries of $P$ lie in $1+v{\R}[v]$.
\end{itemize}
This has the following consequence. Let $\lambda \in \Lambda$
and $1\leq i,j \leq d_\lambda$. For any $h\in \bH_K$, we denote by
$\fX_\lambda^{ij}(h)$ the $(i,j)$-entry of the matrix $\fX_\lambda(h)$.
Then, by \cite[Theorem~4.4 and Remark~4.6]{my02}, we have
\[ v^{a(\lambda)} \fX_\lambda^{ij}(T_w) \in \cO  \qquad \mbox{and}\qquad
v^{a(\lambda)} \fX_\lambda^{ij}(C_w) \in \cO\]
for any $w\in W$ and
\[ v^{a(\lambda)} \fX_\lambda^{ij}(T_w) \equiv v^{a(\lambda)}
\fX_\lambda^{ij}(C_w) \quad \bmod\fp.\]
Hence, the above three elements of $\cO$ have the same constant term which
which we write as $\varepsilon_w\,c_{w,\lambda}^{ij}$. The constants
$c_{w,\lambda}^{ij}\in \R$ are called the {\em leading matrix coefficients}
of $\fX_\lambda$. By \cite[Theorem~4.4]{my02},
these coefficients have the following property:
\begin{alignat*}{2}
c_{w,\lambda}^{ij}&=c_{w^{-1},\lambda}^{ji} &&\qquad
\mbox{for all $w \in W$},\\ c_{w,\lambda}^{ij} &\neq 0 &&\qquad
\mbox{for some $w\in W$}.
\end{alignat*}
We have the following {\em Schur relations}. Let $\lambda,\mu \in \Lambda$, 
$1\leq i,j \leq d_\lambda$ and $1 \leq k,l\leq d_\mu$; then
\[ \sum_{w \in W} c_{w,\lambda}^{ij}c_{w,\mu}^{kl}=\left\{\begin{array}{cl}
\delta_{ik}\,\delta_{jl}\, r_\lambda & \mbox{if $\lambda=\mu$},
\\ 0 & \mbox{if $\lambda\neq \mu$}; \end{array}\right.\]
see \cite[Theorem~4.4]{my02}. (Here $\delta_{ij}$ is the Kronecker delta.) 
The leading matrix coefficients are related to the left cells of $W$ by
the following result. Recall that, given a left cell $\fC$, we have a
corresponding left cell module $[\fC]_A$. Extending scalars from $A$ to
$K$, we obtain an $\bH_K$-module $[\fC]_K$. We denote by $\chi_{\fC}$ the
character of $[\fC]_K$. Then we have:

\addtocounter{thm}{2}
\begin{prop} \label{mylem} Let $\lambda \in \Lambda$ and $\fC$ be a 
left cell of $W$. Then we have 
\[ [\chi_{\fC}:\chi_\lambda] \neq 0 \quad \Leftrightarrow \quad
c_{w,\lambda}^{ij}\neq 0 \mbox{ for some $w\in \fC$ and some $i,j$}.\]
Here, $[\chi_{\fC}:\chi_\lambda]$ denotes the multiplicity of 
$\chi_\lambda$ in $\chi_{\fC}$.
\end{prop}

\begin{proof} By \cite[Prop.~4.7]{my02}, we have the identity 
\[ \sum_{k=1}^{d_\lambda} \sum_{y\in \fC} (c_{y,\lambda}^{ik})^2=
[\chi_{\fC}:\chi_\lambda] \, r_\lambda, \qquad \mbox{for any $1\leq
i \leq d_\lambda$}.\]
Now assume first that $[\chi_{\fC}:\chi_\lambda]\neq 0$. Then, clearly,
we must have $c_{w,\lambda}^{ij}\neq 0$ for some $w\in \fC$ and some
$i,j$. Conversely, assume that $c_{w,\lambda}^{ij}\neq 0$ for some
$w\in \fC$ and some $i,j$. Then the term corresponding to $y=w$ and $k=j$ 
in the sum on the left hand side is non-zero. Since there are no 
cancellations in that sum, the left hand side is non-zero and so
$[\chi_{\fC}:\chi_{\lambda}]\neq 0$, as desired.
\end{proof}

The Schur relations lead to particularly strong results when some additional
hypotheses are satisfied. 

\begin{lem} \label{lem01} Assume that the following condition is satisfied 
for all $\lambda \in \Lambda$:
\begin{equation*}
r_\lambda=1 \quad \mbox{and}\quad c_{w,\lambda}^{ij} \in \Z \quad \mbox{for 
all $w\in W$ and $1\leq i ,j \leq d_\lambda$}.\tag{$*$}
\end{equation*}
Then the following hold.
\begin{itemize}
\item[(a)] We have $c_{w,\lambda}^{ij}\in \{0,\pm 1\}$ for all $w\in W$,
$\lambda \in \Lambda$ and $1\leq i,j \leq d_\lambda$.
\item[(b)] For any $\lambda\in\Lambda$ and $1\leq i,j \leq d_\lambda$,
there exists a unique $w\in W$ such that $c_{w,\lambda}^{ij} \neq 0$; 
we write $w=\bw_\lambda(i,j)$. The correspondence $(\lambda,i,j)\mapsto 
\bw_\lambda(i,j)$ defines a bijective map
\[ \{(\lambda,i,j)\} \mid \lambda \in \Lambda,1\leq i,j \leq d_\lambda\}
\longrightarrow W.\]
\end{itemize}
\end{lem}

\begin{proof} This is proved in \cite[Lemma~3.8]{GeIa05}; see also
\cite[Theorem~4.10]{my02}. In order to illustrate the role of the 
conditions in ($*$), we just recall here the proof of the uniqueness 
statement in (b). Let $\lambda\in \Lambda$ and $1\leq i,j \leq d_\lambda$. 
Then consider the Schur relation where $\lambda=\mu$, $i=l$ and $j=k$. 
This yields:
\[ \sum_{w \in W} (c_{w,\lambda}^{ij})^2=r_\lambda=1.\]
Since the leading matrix coefficients are integers, we conclude that
there exists a unique $w\in W$ such that $c_{w,\lambda}^{ij}=\pm 1$ and 
$c_{y,\lambda}^{ij}=0$ for all $y\in W\setminus \{w\}$. Thus, we have a 
map $(\lambda,i,j)\mapsto w=\bw_\lambda(i,j)$ as in (b). The remaining 
statements are proved in [{\em loc.\ cit.}].
\end{proof}

\begin{rem} \label{incells} In the setting of Lemma~\ref{lem01}, let
$\lambda \in \Lambda$ and set 
\begin{align*}
\fR(\lambda)&:=\{w\in W\mid c_{w,\lambda}^{ij} \neq 0
\mbox{ for some $1\leq i,j \leq d_\lambda$\}}\\
&\;=\{\bw_\lambda(i,j) \mid 1\leq i,j \leq d_\lambda\}.
\end{align*}
It follows from \cite[Theorem~4.4(b)]{my02} that
\begin{equation*}
\bw_\lambda(i,j)^{-1}=\bw_\lambda(j,i).\tag{a}
\end{equation*}
In particular, $\bw_\lambda(i,j)$ is an involution if and only of $i=j$.
Furthermore, let $i,j \in \{1,\ldots,d_\lambda\}$ and define
\begin{equation*}
{^j}\cC_\lambda:=\{\bw_\lambda(k,j)\mid 1\leq k \leq d_\lambda\} \quad
\mbox{and}\quad \cC_\lambda^i:=\{\bw_\lambda(i,l)\mid 1\leq l\leq d_\lambda\}.
\tag{b}
\end{equation*}
It is shown in \cite[Theorem~4.10]{my02} that ${^j}\cC_\lambda$ is contained
in a left cell of $W$ and $\cC_\lambda^j$ is contained in a right cell of $W$.
In particular, the whole set $\fR(\lambda)$ is contained in a two-sided cell
of $W$.
\end{rem}

The following two results will be useful for the identification of 
irreducible characters; they rely on Proposition~\ref{mylem} and 
Lemma~\ref{cwi}.

\begin{lem} \label{indeps} Assume that condition ($*$) in Lemma~\ref{lem01}
is satisfied. Let $I \subseteq S$ and $\lambda \in \Lambda$ be such that
\[ [\Ind_I^S(\varepsilon_I):\chi_\lambda]\neq 0.\]
Then the following hold.
\begin{itemize}
\item[(a)] We have $w \leq_{\cLR} w_I$ for any $w \in \fR(\lambda)$.
\item[(b)] Suppose that $a(\lambda)=l(w_I)$. Then 
\[ w_I \in \fR(\lambda) \qquad \mbox{and}\qquad [\chi_{\fC}:\chi_\lambda]
\neq 0,\]
where $\fC$ is the left cell of $W$ containing $w_I$.
\end{itemize}
\end{lem}

\begin{proof} By Lemma~\ref{cwi}, $X_Iw_I$ is a union of left cells of $W$; 
furthermore, the left cell module $[X_Iw_I]_A$ affords the character
$\Ind_I^S(\varepsilon_I)$. Hence, our assumption that $\chi_\lambda$ occurs 
with non-zero multiplicity in that induced character implies that there exists 
some left cell $\fC$ of $W$ such that
\[ \fC \subseteq X_Iw_I \qquad \mbox{and}\qquad [\chi_{\fC}:\chi_\lambda]
\neq 0.\]
By Proposition~\ref{mylem}, we have $c_{z,\lambda}^{ij}\neq 0$ for some
$z\in \fC$ and some $i,j$. In particular, we have $z \in \fR(\lambda)$. 
Furthermore, writing $z=xw_I$ where $x\in X_I$, we see that $z \leq_{\cL} 
w_I$. Since all elements of $\fR(\lambda)$ are contained in a two-sided 
cell, we obtain (a). 

Now suppose that $a(\lambda)=l(w_I)$. Let $e_I\in \bH_{I,K}$ be the primitive
idempotent affording the sign representation $\varepsilon_I$ of $\bH_{I,K}$. 
By Frobenius reciprocity, we know that $\varepsilon_I$ occurs with non-zero 
multiplicity in the restriction of $\chi_\lambda$ to $\bH_I$. Hence we 
conclude that 
\[ \chi_\lambda(e_I) \in \N \qquad \mbox{and}\qquad \chi_\lambda(e_I)
\neq 0.\]
Now, by Lemma~\ref{cwi}, we have $P_I\, e_I=\varepsilon_{w_I}\, v^{l(w_I)}\, 
C_{w_I}$.  Since $P_I \in 1+\fp$, we obtain 
\begin{align*}
\sum_{i=1}^{d_\lambda} v^{a(\lambda)}\,&\fX_\lambda^{ii}(C_{w_I})\equiv 
v^{l(w_I)}\sum_{i=1}^{d_\lambda} \fX_\lambda^{ii}(C_{w_I})\\&\equiv 
\pm \sum_{i=1}^{d_\lambda} \fX_\lambda^{ii}(e_I)\equiv 
\pm \chi_{\lambda}(e_I)\not\equiv  0 \quad \bmod \fp.
\end{align*}
Hence, there is some $i$ such that $c_{w_I,\lambda}^{ii}\neq 0$ and
so $w_I\in \fR(\lambda)$. Then Proposition~\ref{mylem} shows that
$[\chi_{\fC}:\chi_{\lambda}]\neq 0$, where $\fC$ is the left cell 
containing $w_I$.
\end{proof}

\begin{lem} \label{indeps2} Assume that condition ($*$) in Lemma~\ref{lem01}
is satisfied. Let $I\subseteq S$, $x\in X_I$ and $\lambda \in \Lambda$ be
such that $xw_I \in \fR(\lambda)$. Then we have 
\[ [\Ind_I^S(\varepsilon_I):\chi_\lambda]\neq 0.\]
\end{lem}

\begin{proof} We have $xw_\lambda=\bw_\lambda(i,j)$ where $1\leq i,j
\leq d_\lambda$. Let $\fC$ be the left cell containing $xw_\lambda$.
Then Proposition~\ref{mylem} shows that $[\chi_{\fC}:\chi_{\lambda}]
\neq 0$. On the other hand, $\chi_{\fC}$ is a summand in 
$\Ind_I^S(\varepsilon_I)$, by Lemma~\ref{cwi}. Consequently, $\chi_\lambda$ 
occurs with non-zero-multiplicity in that induced character, as claimed.
\end{proof}

\begin{exmp} \label{symgroup} Let $W=\fS_n$ be the symmetric group on 
$\{1,\ldots,n\}$. Then $W$ is a Coxeter group with generating set
$S=\{s_1,\ldots,s_{n-1}\}$, where $s_i=(i,i+1)$ for $1\leq i \leq n-1$.  The
diagram is given as follows.
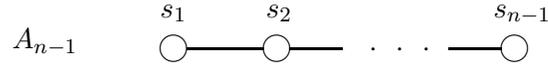
\begin{center}
\begin{picture}(300,30)
\put( 40, 5){$A_{n-1}$}
\put(101, 5){\circle{10}}
\put(106, 5){\line(1,0){29}}
\put(140, 5){\circle{10}}
\put(145, 5){\line(1,0){20}}
\put(175, 2){$\cdot$}
\put(185, 2){$\cdot$}
\put(195, 2){$\cdot$}
\put(205, 5){\line(1,0){20}} \put(230, 5){\circle{10}}
\put( 96, 17){$s_1$}
\put(136, 17){$s_2$}
\put(222, 17){$s_{n-1}$}
\end{picture}
\end{center}
Let $\bH=\bH_A(\fS_n,S)$ be the corresponding Iwahori--Hecke algebra. Let 
$\Lambda_n$ be the set of all partitions of $n$. Then it is well-known 
(see, for example, Hoefsmit \cite{Hoefs} or Dipper--James \cite{DiJa1}) 
that we have a labelling 
\[\Irr(\bH_K)=\{\chi_{\lambda} \mid \lambda \in \Lambda_n\}.\]
It will be important to specify precisely the labelling that we are using. 
A ``standard'' labelling is defined in \cite[\S 5.4]{ourbuch}; however, in 
order to avoid going back and forth between partitions and their conjugates,
we shall change that labelling and simply denote here by $\chi_\lambda$ 
the irreducible character which is labelled by the conjugate partition 
$\lambda^*$ in [{\em loc.\ cit.}]. (See also Corollary~\ref{ident} where 
we describe a representation affording $\chi_\lambda$.) For example, with 
this convention, $\chi_{(1^n)}$ is the trivial character and $\chi_{(n)}$
is the sign character of $\bH_K$.  The labelling that we have chosen is
characterized as follows in terms of induced characters (and this
characterisation will be important in Section~4).

Let $\lambda \in \Lambda_n$. If $\lambda$ has parts $\lambda_1\geq 
\lambda_2 \geq \cdots \geq 0$, we set 
\[ \lambda_i^+=\lambda_1+\lambda_2+\cdots +\lambda_i \qquad \mbox{for
$i=1,2,3,\ldots$}.\]
Let $I_\lambda:=\{1,\ldots,n\}\setminus \{\lambda_i^+\mid i=1,2,\ldots\}$. 
Then $I_\lambda\subseteq S$ and $\fS_\lambda:=W_{I_\lambda}$ is the 
Young subgroup corresponding to $\lambda$. Let $\varepsilon_\lambda$ be 
the sign character of the corresponding parabolic subalgebra $\bH_\lambda 
\subseteq \bH$. Then, by \cite[5.4.7]{ourbuch}, the labelling of 
$\Irr(\bH_K)$ is uniquely determined by the condition that
\begin{equation*}
\Ind_{I_\lambda}^{S}(\varepsilon_\lambda)=\chi_{\lambda}+ \mbox{sum of 
characters $\chi_\mu$ where $\lambda \triangleleft \mu$},\tag{a}
\end{equation*}
where $\trianglelefteq$ denotes the usual dominance order on $\Lambda_n$.
(We have $\lambda \trianglelefteq \mu$ if and only if $\lambda_i^+\leq\mu_i^+$
for $i=1,2,3,\ldots$; we write $\lambda \triangleleft \mu$ if $\lambda
\trianglelefteq \mu$ and $\lambda\neq \mu$.) 

The formula for the Schur elements in \cite[Prop.~9.4.5]{ourbuch} and the
identities in \cite[\S 5.4]{ourbuch} now show that we have 
\begin{equation*}
r_\lambda=1 \qquad \mbox{and}\qquad a(\lambda)=l(w_\lambda) \qquad 
\mbox{for any $\lambda \in \Lambda_n$},\tag{b}
\end{equation*}
where $w_\lambda$ is the unique element of maximal length in $\fS_\lambda$. 
Furthermore, by the discussion in \cite[\S 5]{my02}, each $\chi_\lambda$ is 
afforded by an orthogonal representation such that the leading matrix 
coefficients $c_{w,\lambda}^{ij}$ are integers. (Since that discussion 
is somewhat sketchy, a more rigorous argument is given in 
\cite[Remark~3.7]{GeIa05}, based on the Dipper--James construction 
of Hoefsmit's matrices in \cite[Theorem~4.9]{DiJa2}.) 
Thus, condition ($*$) in Lemma~\ref{lem01} is satisfied. Hence, each 
element of $\fS_n$ lies in $\fR(\lambda)$ for a unique $\lambda \in 
\Lambda_n$. Furthermore, applying Lemma~\ref{indeps} and using (a) and (b), 
we obtain that $w_\lambda \in \fR(\lambda)$.  Thus, we have a partition 
\begin{equation*}
\fS_n=\coprod_{\lambda\in\Lambda_n} \fR(\lambda) \qquad \mbox{where}
\qquad w_\lambda \in \fR(\lambda).\tag{c}
\end{equation*}
\end{exmp}

\begin{rem} \label{choice} In the setting of Example~\ref{symgroup},
we will henceforth fix, for each $\lambda \in \Lambda_n$, an orthogonal 
representation $\fX_\lambda \colon \bH_K \rightarrow M_{d_\lambda}(K)$ 
affording $\chi_\lambda$ such that the leading matrix coefficients 
$c_{w,\lambda}^{ij}$ are integers. Thus, the sets $\fR(\lambda)$ are 
defined and we have $w_\lambda \in \fR(\lambda)$. By Remark~\ref{incells}, 
we have 
\[ \fR(\lambda)=\{\bw_\lambda(i,j) \mid 1\leq i,j\leq d_\lambda\}.\]
Since $w_\lambda\in \fR(\lambda)$ is an involution, we have
$w_\lambda=\bw_\lambda(i,i)$ for some $i$. Conjugating $\fX_\lambda$ by
a suitable permutation matrix, we may in fact assume that $i=1$. In the 
sequel, we shall always assume that $\fX_\lambda$ has been ``normalized''
in this way, that is, we have 
\begin{equation*}
w_\lambda=\mathbf{w}_\lambda(1,1) \in \fR(\lambda).  \tag{a}
\end{equation*}
Let $X_\lambda$ be the set of distinguished left coset representatives 
of $\fS_\lambda$ in $\fS_n$. Now, by Remark~\ref{incells}, the set
${^1}\cC_\lambda=\{\bw_\lambda(i,1)\mid 1\leq i \leq d_\lambda\}$ is 
contained in a left cell of $\fS_n$. Hence, by Lemma~\ref{cwi}, we have 
\begin{equation*}
\bw_\lambda(i,1)=x_i w_\lambda\qquad\mbox{where $x_i\in X_\lambda$ for
$1 \leq i \leq d_\lambda$ and $x_1=1$}.\tag{b}
\end{equation*}
(An explicit description of the elements $x_i$ will be given in 
Remark~\ref{palli}.) 

We shall need the following result.
\end{rem}

\begin{lem} \label{technic} Let $\lambda \in \Lambda_n$ and $x\in
X_\lambda$. Assume that $\mu \in \Lambda_n$ is such that $xw_\lambda 
\in \fR(\mu)$. Then we have $\lambda \trianglelefteq \mu$; furthermore,
if $\mu=\lambda$, then $x=x_i$ for some $i\in \{1,\ldots,d_\lambda\}$.
\end{lem}

\begin{proof} By Lemma~\ref{indeps2}, $\chi_\mu$ occurs with non-zero 
multiplicity in $\Ind_{I_\lambda}^S(\varepsilon_\lambda)$. Hence 
Example~\ref{symgroup}(a) shows that $\lambda \trianglelefteq \mu$.
Thus, it remains to consider the case where $\lambda=\mu$. Now, by
Lemma~\ref{cwi}, the set $X_\lambda w_\lambda$ is a union of left cells
of $\fS_n$. Let $\fC,\fC'\subseteq X_\lambda w_\lambda$ be left cells such
that $w_\lambda \in \fC$ and $xw_\lambda \in \fC'$. First we claim that
$\fC=\fC'$. This is seen as follows. By Lemma~\ref{indeps}(b), 
$\chi_\lambda$ occurs with non-zero multiplicity in $\chi_{\fC}$.
On the other hand, since $xw_\lambda \in \fR(\lambda)$, we can apply 
Proposition~\ref{mylem} to conclude that $\chi_\lambda$ also occurs
with non-zero multiplicity in $\chi_{\fC'}$. Now assume, if possible,
that $\fC \neq \fC'$. Since $\chi_{\fC}$ and $\chi_{\fC'}$ both occur 
as summands in $\Ind_{I_\lambda}^S(\varepsilon_\lambda)$ (see
Lemma~\ref{cwi}), we conclude that $\chi_\lambda$ occurs with 
multiplicity at least $2$ in that induced character, contradicting
Example~\ref{symgroup}(a). Thus, we must have $\fC=\fC'$, as claimed.

Now let us write $xw_\lambda=\bw_\lambda(i,j)\in \fC$ where $1\leq i,j
\leq d_\lambda$. As in the proof of Proposition~\ref{mylem}, we have 
the identity 
\[ \sum_{k=1}^{d_\lambda} \sum_{y \in \fC} (c_{y,\lambda}^{ik})^2=
[\chi_{\fC}:\chi_\lambda]\, r_\lambda=1.\]
As discussed in Remark~\ref{choice}, we have 
$x_iw_\lambda=\bw_\lambda(i,1) \in \fC$.  Hence this term gives a
non-zero contribution to the sum on the left hand side. On the
other hand, we also have $c_{xw_\lambda,\lambda}^{ij}\neq 0$ and
this gives a non-zero contribution. Since the right hand side equals $1$, 
we conclude that $xw_\lambda=x_iw_\lambda$ and so $x=x_i$, as desired.
\end{proof}

\begin{rem} \label{kostka} Let $\lambda,\mu \in \Lambda_n$ be such
that $\lambda \trianglelefteq \mu$. Then the multiplicity  
\[ [\Ind_{I_\lambda}^S(\varepsilon_\lambda):\chi_\mu]\]
can be expressed in a purely combinatorial way, by {\em Young's rule}.
In particular, it is known that the above multiplicity is non-zero; see, for
example, Murphy \cite[Theorem~7.2]{Mu2}. So Lemma~\ref{indeps} shows 
that $w \leq_{\cLR} w_\lambda$ for any $w \in \fR(\mu)$. Using 
Example~\ref{symgroup}(c), we conclude that 
\begin{equation*}
\lambda \trianglelefteq \mu \qquad \Rightarrow \qquad w_\mu
\leq_{\cLR} w_\lambda.\tag{a}
\end{equation*}
In Corollary~\ref{ourmain1}, we will see that the converse also holds, but 
the proof requires much more work. 
\end{rem}

\section{The Murphy basis} \label{sec-mb}

Throughout this and the following section, we consider the Iwahori--Hecke 
algebra $\bH=\bH_A(\fS_n,S)$, as in Example~\ref{symgroup}. The conventions
in Remark~\ref{choice} will also remain in force. 

In two fundamental articles \cite{Mu1} and \cite{Mu2}, Murphy  has 
constructed a new basis $\{x_{\fs\ft}\}$ of $\bH$ whose elements are 
indexed by pairs of standard $\lambda$-tableaux, for various partitions 
$\lambda$ of $n$. As we shall see, the elements $x_{\fs\ft}$ are defined 
as a mixture of basis elements $T_w$ and $C_w$. Our aim is to establish a 
direct link with the Kazhdan--Lusztig basis; this will be achieved in 
Theorem~\ref{ourmain}.

We begin by recalling some purely combinatorial notions, where we follow 
\cite{Mu2} (but we let $\fS_n$ act on the left on $\{1,\ldots,n\}$). Another
reference is the exposition by Mathas \cite{Mathas99}.

As before, let $\Lambda_n$ be the set of all partitions of $n$. Let $\lambda
\in \Lambda_n$ and let $\lambda_1 \geq \ldots \geq \lambda_r>0$ be the 
non-zero parts of $\lambda$. The correspdonding Young diagram $[\lambda]$ is 
the set of all pairs $(i,j)$ such that $1\leq i \leq r$ and $1\leq j \leq 
\lambda_i$. A $\lambda$-tableau is a bijection $\ft\colon [\lambda]
\rightarrow \{1,\ldots,n\}$. We say that $\ft$ is {\em row-standard} if the 
sequence $\ft(i,1),\ft(i,2),\ldots$ is strictly increasing for each $i$; 
similarly, we say that $\ft$ is {\em column-standard} if the sequence
$\ft(1,j),\ft(2,j),\ldots$ is strictly increasing for each~$j$. We say 
that $\ft$ is {\em standard} if $\ft$ is row-standard and column-standard.
We denote by $\ft^\lambda$ the unique standard $\lambda$-tableau such that 
\[ \ft^\lambda(i,j)=\lambda_1+\lambda_2+\cdots +\lambda_{i-1}+j\quad
\mbox{for $1\leq i\leq r$ and $1\leq j \leq \lambda_i$}.\]
Here are some examples, where $n=5$ and $\lambda=(3,2)$:
\begin{center}
\begin{picture}(280,50)
\put(00,27){$\ft^\lambda=$}
\put(30,15){\line(1,0){30}}
\put(30,30){\line(1,0){45}}
\put(30,45){\line(1,0){45}}
\put(30,15){\line(0,1){30}}
\put(45,15){\line(0,1){30}}
\put(60,15){\line(0,1){30}}
\put(75,30){\line(0,1){15}}
\put(35,34){$1$}
\put(50,34){$2$}
\put(65,34){$3$}
\put(35,19){$4$}
\put(50,19){$5$}
\put(130,0){$\text{\small (standard)}$}
\put(130,15){\line(1,0){30}}
\put(130,30){\line(1,0){45}}
\put(130,45){\line(1,0){45}}
\put(130,15){\line(0,1){30}}
\put(145,15){\line(0,1){30}}
\put(160,15){\line(0,1){30}}
\put(175,30){\line(0,1){15}}
\put(135,34){$1$}
\put(150,34){$3$}
\put(165,34){$5$}
\put(135,19){$2$}
\put(150,19){$4$}
\put(230,0){$\text{\small (row-standard)}$}
\put(230,15){\line(1,0){30}}
\put(230,30){\line(1,0){45}}
\put(230,45){\line(1,0){45}}
\put(230,15){\line(0,1){30}}
\put(245,15){\line(0,1){30}}
\put(260,15){\line(0,1){30}}
\put(275,30){\line(0,1){15}}
\put(235,34){$1$}
\put(250,34){$4$}
\put(265,34){$5$}
\put(235,19){$2$}
\put(250,19){$3$}
\end{picture}
\end{center}
The group $\fS_n$ acts naturally on $\lambda$-tableaux, the action being 
given by $(w.\ft)(i,j):=w.\ft(i,j)$ for a tableau $\ft$ and $w\in \fS_n$. 
The row stabiliser of $\ft^\lambda$ is the Young subgroup $\fS_\lambda$.
As in Remark~\ref{choice}, let $X_\lambda$ be the set of distinguished 
left coset representatives of $\fS_\lambda$ in $\fS_n$. By Dipper--James 
\cite[Lemma~1.1]{DiJa1}, we have the following explicit description:
\[X_\lambda=\{w\in\fS_n\mid w.\ft^\lambda \mbox{ is row-standard}\}.\]
As in [{\em loc. cit.}], if $\ft$ is a row-standard $\lambda$-tableau, 
the unique element $d\in X_\lambda$ such that $\ft=d.\ft^\lambda$ will 
be denoted by $d(\ft)$.

Let $\bT(\lambda)$ be the set of standard $\lambda$-tableaux.

\begin{thm}[Murphy \protect{\cite{Mu1}, \cite{Mu2}}] \label{murphy12}
For any $\lambda \in \Lambda_n$ and $\fs, \ft\in \bT(\lambda)$, we define 
elements of $\bH$ by 
\[ x_\lambda:=\sum_{w\in \fS_\lambda} v^{l(w)}\,T_w \qquad \mbox{and}
\qquad x_{\fs\ft}:=T_{d(\fs)}\,x_\lambda \,T_{d(\ft)^{-1}}.\]
Then the following hold.
\begin{itemize}
\item[(a)] The set $\{ x_{\fs\ft} \mid \mbox{$\fs,\ft \in \bT(\lambda)$ 
for some $\lambda \in \Lambda_n$}\}$ is an $A$-basis 
of $\bH$.
\item[(b)] For any $\lambda\in \Lambda_n$, let $N^\lambda$ be the 
$A$-submodule of $\bH$ spanned by all elements $x_{\fs\ft}$ where $\fs,\ft
\in \bT(\mu)$ for some $\mu\in \Lambda$ such that $\lambda \trianglelefteq 
\mu$. Then $N^\lambda$ is a two-sided ideal in $\bH$.
\end{itemize}
\end{thm}

The statement in (a) can be found in \cite[Theorem~3.9]{Mu1} or, with a 
somewhat different proof, in \cite[Theorem~4.17]{Mu2}. The statement in (b) 
is proved in \cite[Theorem~4.18]{Mu2} (see also \cite[\S 5]{Mu2}). Note that 
the element that we denote by $T_w$ corresponds to the element $v^{-l(w)}T_w$
in Murphy's notation. Thus, the element denoted by $x_\lambda$ in the above 
statement is exactly the same as in Murphy's work; the element denoted by 
$x_{\fs\ft}$ is only the same up to a power of~$v$. However, this does not 
affect the validity of (a) and (b) since $v$ is invertible in $A$. 

Murphy also obtains the following result concerning the Specht modules of 
$\bH$. For any $\lambda \in \Lambda_n$, let $\hat{N}^\lambda$ be the 
$A$-submodule of $\bH$ spanned by all $x_{\fs\ft}$ where $\fs,\ft\in 
\bT(\mu)$ for some $\mu\in \Lambda$ such that $\lambda \triangleleft \mu$, 
that is, $\lambda \trianglelefteq \mu$ and $\lambda\neq\mu$. Thus, we have 
\[ \hat{N}^\lambda=\sum_{\mu} N^\mu\]
where the sum runs over all $\mu\in \Lambda_n$ such that $\lambda 
\triangleleft \mu$. In particular, $\hat{N}^\lambda$ is a two-sided ideal 
and we have $N^\lambda=\bH x_\lambda\bH +\hat{N}^\lambda$.

\begin{thm}[Murphy \protect{\cite[\S 5]{Mu2}}] \label{murphiA} Let $\lambda 
\in \Lambda_n$. Then 
\[ \tilde{\cS}^\lambda:=\langle T_{d(\fs)}\,x_\lambda+\hat{N}^\lambda 
\mid \fs \in \bT(\lambda)\rangle_A\subseteq N^\lambda/\hat{N}^\lambda\]
is a left $\bH$-module, and $N^\lambda/\hat{N}^\lambda$ is a direct sum
of $|\bT(\lambda)|$ copies of $\tilde{\cS}^\lambda$. Furthermore, the 
$\bH_K$-module $K\otimes_A \tilde{\cS}^\lambda$ is simple, and 
$\tilde{\cS}^\lambda$ is the contragredient dual of the Dipper--James
Specht module defined in \cite{DiJa1}.
\end{thm}

The following reformulation of Theorem~\ref{murphy12} will eventually 
allow us to establish a direct connection between the Murphy basis and 
the Kazhdan--Lusztig  basis.

\begin{cor} \label{murphib} For any $\lambda\in \Lambda_n$ and $\fs,\ft
\in \bT(\lambda)$, we set 
\[ \ty_{\fs\ft}:=T_{d(\fs)}\, C_{w_\lambda} \, T_{d(\ft)^{-1}} \in \bH.\]
Then $\ty_{\fs\ft}=\pm v^{l(w_\lambda)}j(x_{\fs\ft})$ for all 
$\fs,\ft$. Consequently, the following hold.
\begin{itemize}
\item[(a)] The set $\{ \ty_{\fs\ft} \mid \mbox{$\fs,\ft\in \bT(\lambda)$ 
for some $\lambda \in \Lambda_n$}\}$ is an $A$-basis of $\bH$.
\item[(b)] For any $\lambda\in \Lambda_n$, let $\fN^\lambda$ be the 
$A$-submodule of $\bH$ spanned by all elements $\ty_{\fs\ft}$ where $\fs,
\ft\in \bT(\mu)$ for some $\mu\in \Lambda$ such that $\lambda
\trianglelefteq \mu$. Then $\fN^\lambda$ is a two-sided ideal in $\bH$.
\end{itemize}
\end{cor}

\begin{proof} Using the expression for $C_{w_\lambda}'$ in the proof of 
Lemma~\ref{cwi}, we immediately see that $x_\lambda=v^{l(w_\lambda)}\,
C_{w_\lambda}'$. Next recall that, for any $w\in \fS_n$, we have $j(T_w)=
\varepsilon_w T_w$ and $j(C_w')=\varepsilon_w C_w$. Hence we obtain
\begin{align*}
j(x_{\fs\ft})&=j(T_{d(\fs)})\,j(x_\lambda)\,j(T_{d(\ft)^{-1}})\\
&=\varepsilon_{d(\fs)}\, \varepsilon_{d(\ft)}\, \varepsilon_{w_\lambda}\,
v^{-l(w_\lambda)}T_{d(\fs)}C_{w_\lambda}T_{d(\ft)}\\
&=\pm v^{-l(w_\lambda)}\ty_{\fs\ft}.
\end{align*}
Since $j \colon \bH \rightarrow \bH$ is a ring involution, we now see
that the statements in Theorem~\ref{murphy12} hold with $x_{\fs\ft}$ 
replaced by $\ty_{\fs\ft}$ throughout. 
\end{proof}

\begin{rem} \label{ideal} For any $\lambda \in \Lambda_n$, let 
$\hat{\fN}^\lambda$ be the $A$-submodule of $\bH$ spanned by all elements 
$\ty_{\fs\ft}$ where $\fs,\ft\in \bT(\mu)$ for some  $\mu\in \Lambda$ such 
that $\lambda \triangleleft \mu$.  Thus, we have 
\[ \hat{\fN}^\lambda=\sum_{\mu} \fN^\mu\]
where the sum runs over all $\mu\in \Lambda_n$ such that $\lambda 
\triangleleft \mu$. In particular, $\hat{\fN}^\lambda$ is a two-sided ideal 
and we have $\fN^\lambda=\bH C_{w_\lambda}\bH +\hat{\fN}^\lambda$.
We claim that 
\[ j(N^\lambda)=(N^\lambda)^\dagger=\fN^\lambda \qquad \mbox{and}\qquad 
j(\hat{N}^\lambda)=(\hat{N}^\lambda)^\dagger=\hat{\fN}^\lambda.\]
Indeed, in the proof of Corollary~\ref{murphib}, we have seen that
$j(x_\lambda)=\pm v^{-l(w_\lambda)} C_{w_\lambda}$. But then we also have 
\[ x_\lambda^\dagger=\overline{j(x_\lambda)}=\pm v^{l(w_\lambda)}
\overline{C}_{w_\lambda}=\pm v^{l(w_\lambda)}\, C_{w_\lambda}.\]
An easy induction on the dominance order, using the relations
\[ N^\lambda=\bH x_\lambda \bH +\hat{N}^\lambda \qquad \mbox{and}\qquad
\fN^\lambda=\bH C_{w_\lambda} \bH +\hat{\fN}^\lambda,\]
now shows that $(N^\lambda)^\dagger=\fN^\lambda=j(N^\lambda)$. 
Since this holds for all $\lambda \in \Lambda_n$, we also have 
$(\hat{N}^\lambda)^\dagger=\hat{\fN}^\lambda=j(\hat{N}^\lambda)$. 
Thus, the above claim is proved. 
\end{rem}

The identification of $\chi_\lambda$ in the following result relies on 
our conventions for labelling the irreducible characters in 
Example~\ref{symgroup}.

\begin{cor} \label{ident} Let $\lambda \in \Lambda_n$. Then 
\[{\cE}^\lambda:=\langle T_{d(\fs)}\,C_{w_\lambda} +\hat{\fN}^\lambda
\mid\fs\in \bT(\lambda)\rangle_A\subseteq \fN^\lambda/\hat{\fN}^\lambda\]
is a left $\bH$-module, and $\fN^\lambda/\hat{\fN}^\lambda$ is a direct 
sum of $|\bT(\lambda)|$ copies of ${\cE}^\lambda$. The character 
afforded by $K \otimes_A \cE^\lambda$ is $\chi_\lambda$, and we have
$d_\lambda=|\bT(\lambda)|$.
\end{cor}

\begin{proof} By Remark~\ref{ideal}, the ring involution $j\colon \bH
\rightarrow \bH$ transforms $N^\lambda$ and $\hat{N}^\lambda$ into
$\fN^\lambda$ and $\hat{\fN}^\lambda$, respectively. Hence we have an 
induced isomorphism of additive groups 
\[ \theta\colon N^\lambda/\hat{N}^\lambda \rightarrow \fN^\lambda/
\hat{\fN}^\lambda\]
such that $\theta(h.e)=j(h).\theta(e)$ for all $h\in \bH$ and $e\in
N^\lambda/\hat{N}^\lambda$. Now, by Corollary~\ref{murphib}, we have 
$j(x_\lambda)=\pm v^{-l(w_\lambda)}\, C_{w_\lambda}$ and so 
$\theta(\tilde{\cS}^\lambda)=\cE^\lambda$. Hence Theorem~\ref{murphiA} 
implies that $\cE^\lambda$ is a left $\bH$-module such that 
$K \otimes_A \cE^\lambda$ is simple; furthermore, $\fN^\lambda/
\hat{\fN}^\lambda$ is a direct sum of $|\bT(\lambda)|$ copies of 
$\cE^\lambda$. 

Thus, we see that it will be sufficient to prove that the character afforded 
by $(\fN^\lambda)/\hat{\fN}^\lambda)_K$ is $d_\lambda \chi_\lambda$. 
(The subscript $K$ indicates extension of scalars from $A$ to $K$.) 
For this purpose, we proceed by downward induction on the dominance order.
The unique maximal element for that order is the partition $\lambda=(n)$.
In this case, we have $\fS_{(n)}=\fS_n$ and $X_{(n)}=\{1\}$. Hence we
have $\fN^{(n)}=\langle C_{w_{(n)}}\rangle_A$ and $\hat{\fN}^{(n)}=\{0\}$.
By Lemma~\ref{cwi}, $\fN^{(n)}$ affords the sign character. On the
other hand, by Example~\ref{symgroup}, $\chi_{(n)}$ also is the sign
character and $d_{(n)}=1$. Thus, the assertion holds in this case. Now 
let $\lambda\neq (n)$. First we prove the following two statements.
\begin{itemize}
\item[(1)] {\em If $\lambda\triangleleft \mu$, then $\chi_\mu$ does not 
occur in the character of $(\fN^\lambda/\hat{\fN}^\lambda)_K$}.
\item[(2)] {\em Consider the left $\bH$-module $M^\lambda:=
\bH C_{w_\lambda}$. Then $(\fN^\lambda/ \hat{\fN}^\lambda)_K$ and 
$M^\lambda_K$ have a common simple constituent}.
\end{itemize}
To prove (1), we argue as follows. As a left $\bH_K$-module, we have 
\[ \fN^\lambda_K \cong (\fN^\lambda/\hat{\fN}^\lambda)_K \oplus 
\hat{\fN}^\lambda_K \qquad \mbox{where}\qquad \fN_K^\mu\subseteq 
\hat{\fN}^\lambda_K.\]
By induction, we already know that $d_\mu\chi_\mu$ is the character 
afforded by $({\fN}^\mu/\hat{\fN}^\mu)_K$. Hence, the character
of $\hat{\fN}^\lambda$ contains $d_\mu \chi_\mu$ as a summand. Since 
$d_\mu$ is the multiplicity of $\chi_\mu$ in the character of the 
regular representation of $\bH_K$ and $\fN_K^\lambda$ is contained
in $\bH_K$, we conclude that (1) holds. Now let us prove (2). Since 
$C_{w_\lambda} \in \fN^\lambda$ and $C_{w_\lambda} \not\in 
\hat{\fN}^\lambda$, the inclusion $M^\lambda \subseteq \fN^\lambda$ 
induces a non-zero homomorphism of $\bH$-modules $\varphi \colon
M^\lambda \rightarrow \fN^\lambda/\hat{\fN}^\lambda$. Hence, (2) follows. 

Now we can determine the character afforded by $(\fN^\lambda/
\hat{\fN}^\lambda)_K$. Let $\mu \in \Lambda_n$ be such that $\chi_\mu$ 
is the character of a common simple component as in (2). By 
Lemma~\ref{cwi}, we have $M^\lambda \cong \Ind_{I_\lambda}^S
(\varepsilon_\lambda)$. Hence, since $\chi_\mu$ occurs in $M^\lambda_K$, 
we must have $\lambda \trianglelefteq \mu$ by Example~\ref{symgroup}(a). 
But then (1) yields that $\lambda=\mu$. Thus, $\chi_\lambda$ occurs in 
the character of $(\fN^\lambda/\hat{\fN}^\lambda)_K$. Since the latter 
module is a direct sum of $|\bT(\lambda)|$ copies of $\cE_K^\lambda$, we 
conclude that $\chi_\lambda$ is the character of $\cE_K^\lambda$. Since 
$\cE_K^\lambda$ has dimension $|\bT(\lambda)|$, this also yields that 
$|\bT(\lambda)|= d_\lambda$, as required.
\end{proof}

\begin{cor} \label{nonzero} Let $\lambda,\mu \in \Lambda_n$ and $h\in
\fN^\lambda$. Then $\fX_\mu(h)=0$ unless $\lambda \trianglelefteq \mu$.
\end{cor}

\begin{proof} Let us assume that $\fX_\mu(h)\neq 0$. We must show that
$\lambda \trianglelefteq \mu$. To see this, consider the representation
afforded by the left $\bH$-module $\fN^\mu/\hat{\fN}^\mu$. Extending 
scalars to $K$ and using Corollary~\ref{ident}, we see that this 
representation is equivalent to a direct sum of $d_\mu$ copies of $\fX_\mu$.
So the condition $\fX_\mu(h)\neq 0$ implies that $h.(\fN^\mu/\hat{\fN}^\mu)
\neq 0$. This means that there exist some standard $\mu$-tableaux $\fs,
\fs_1,\ft, \ft_1$ such that $\ty_{\fs_1 \ft_1}$ occurs with non-zero 
coefficient in the decomposition of $h \ty_{\fs\ft}$. Now, since $h\in 
\fN^\lambda$ and $\fN^\lambda$ is a two-sided ideal, we can write 
$h\ty_{\fs\ft}$ as a linear combination of terms $\ty_{\fu\fv}$ where 
$\fu, \fv$ are standard $\nu$-tableaux for partitions $\nu \in\Lambda_n$
such that $\lambda \trianglelefteq\nu$. One of the terms in that linear 
combination is $\ty_{\fs_1\ft_1}$. Hence, we conclude that $\lambda
\trianglelefteq \mu$, as claimed.
\end{proof}

The following definition relies on our conventions in Remark~\ref{choice} 
concerning the labelling of the elements in the sets $\fR(\lambda)$.

\begin{defn} \label{maindef} Let $w\in\fS_n$. Let $\lambda \in \Lambda_n$ be
such that $w\in\fR(\lambda)$. Then $w=\bw_\lambda(i,j)$ where $1\leq i,j\leq 
d_\lambda$ are uniquely determined. We set
\[ \lambda_w:=\lambda,\qquad P_\lambda:=P_{I_\lambda},\qquad 
\alpha_w:=a(\lambda)=l(w_\lambda)\]
(see Example~\ref{symgroup}) and 
\[Z_w:=\frac{1}{P_\lambda}\, \varepsilon_{w_\lambda}\,
v^{l(w_\lambda)}\, C_{x_iw_\lambda} \,C_{w_\lambda x_j^{-1}}\in \bH_K.\]
Note that we do need to extend scalars from $A$ to $K$ in order to 
define the above element. 
\end{defn}

The following result shows that $Z_w$ actually lies in $\bH$.

\begin{lem} \label{mainlem1} We have $Z_w\in \fN^\lambda\subseteq \bH$
and $\overline{Z}_w= Z_w$. For $w=w_\lambda$, we have
$Z_{w_\lambda}=C_{w_\lambda}$.
\end{lem}

\begin{proof} By \cite[Prop.~3.3]{myind}, we have 
\begin{align*}
C_{x_iw_\lambda}&=\mbox{ $A$-linear combination of terms $T_xC_{w_\lambda}$ 
where $x\in X_\lambda$},\\
C_{w_\lambda x_j^{-1}}&=\mbox{ $A$-linear combination of terms $C_{w_\lambda}
T_{y^{-1}}$ where $y\in X_\lambda$}.
\end{align*}
(We have already used this fact in the proof of Lemma~\ref{cwi}.) Hence, 
$Z_w$ is an $A$-linear combination of terms of the form 
\[\frac{1}{P_\lambda}\, \varepsilon_{w_\lambda}\, v^{l(w_\lambda)}\, 
T_x \,C_{w_\lambda}^2\, T_{y^{-1}}\qquad\mbox{where $x,y\in X_\lambda$}.\]
Using Lemma~\ref{cwi}, each of the above terms is equal to 
$T_xC_{w_\lambda}T_{y^{-1}}$ and, hence, lies in $\bH$. Consequently,
we also have that $Z_w$ lies in $\bH C_{w_\lambda} \bH$, and this is
contained in $\fN^\lambda$ by Remark~\ref{ideal}. Finally, note that 
$v^{l(w_\lambda)}\overline{P}_\lambda =v^{-l(w_\lambda)} P_\lambda$. This 
yields $\overline{Z}_w=Z_w$. 
\end{proof}

\begin{lem} \label{mainlem2} Let $w\in \fR(\lambda)$ and 
$1\leq k,l\leq d_\lambda$. Then we have
\[ v^{a(\lambda)}\,\fX_\lambda^{kl}(Z_w) \in \cO \qquad \mbox{and}\qquad
v^{a(\lambda)}\,\fX_\lambda^{kl}(Z_w) \equiv \pm \delta_{ik}\delta_{jl}
\quad \bmod \fp,\]
where $i,j \in \{1,\ldots,d_\lambda\}$ are such that $w=\bw_\lambda(i,j)$.
\end{lem}

\begin{proof} We write $w=\bw_\lambda(i,j)$ where $1\leq i,j\leq d_\lambda$.
We have $a(\lambda)=l(w_\lambda)$  by Example~\ref{symgroup}. 
Hence we obtain
\begin{align*}
v^{a(\lambda)}\,\fX_\lambda^{kl}(Z_w) &=\pm \frac{1}{P_\lambda}\,
v^{2a(\lambda)}\, \fX_{\lambda}^{kl}\Bigl(C_{x_iw_\lambda}\, 
C_{w_\lambda x_j^{-1}}\Bigr) \\
&=\pm \frac{1}{P_\lambda}\,\sum_{r=1}^{d_\lambda} \Bigl(v^{a(\lambda)}\, 
\fX_{\lambda}^{kr}(C_{x_iw_\lambda})\Bigr)\, \Bigl(v^{a(\lambda)}\, 
\fX_\lambda^{rl}(C_{w_\lambda x_j^{-1}})\Bigr).
\end{align*}
First of all, this shows that the above expression lies in $\cO$; note
that $P_\lambda \in 1+v{\Z}[v]$. Furthermore, its constant term can be 
expressed by the leading matrix coefficients of $x_iw_\lambda$ and 
$x_jw_\lambda$. Indeed, since $x_i w_\lambda=\bw_\lambda(i,1)$,
Lemma~\ref{lem01} shows that 
\[v^{a(\lambda)}\,\fX_{\lambda}^{kr}(C_{x_iw_\lambda})\equiv 
\varepsilon_{x_iw_\lambda}\, c_{x_iw_\lambda,\lambda}^{kr} \equiv \pm 
\delta_{ki} \delta_{r1} \quad\bmod \fp.\]
Similarly, since $w_\lambda x_j^{-1}=\bw_\lambda(1,j)$, we have 
\[v^{a(\lambda)}\,\fX_{\lambda}^{rl}(C_{w_\lambda x_j^{-1}})\equiv\pm 
\delta_{lj} \delta_{r1} \quad\bmod \fp.\]
Now note that $P_\lambda \in 1+\fp$ and so $P_\lambda^{-1} \equiv 1 
\bmod \fp$.   Hence, we obtain 
\[ v^{a(\lambda)}\,\fX_\lambda^{kl}(Z_w) \equiv \pm \delta_{ki}
\delta_{lj} \quad \bmod\fp,\]
as desired.
\end{proof}

\begin{thm} \label{ourmain} Let $\lambda\in\Lambda_n$. Then the following
hold.
\begin{align*}
\eta_w\, C_w&\in Z_w +\hat{\fN}^\lambda\subseteq \fN^\lambda \quad 
\mbox{for all $w\in \fR(\lambda)$, where $\eta_w=\pm 1$},\tag{a}\\
\fN^\lambda& =\langle C_w \mid w\in \fR(\mu) \mbox{ for some $\mu\in 
\Lambda_n$ such that $\lambda \trianglelefteq \mu$}\rangle_A.\tag{b}
\end{align*}
\end{thm}

(The sign in (a) will be determined explicitly in Corollary~\ref{sign}.)

\begin{proof} (a) We proceed by downward induction on the dominance order. 
The unique maximal element for that order is the partition $\lambda=(n)$.
In this case, we have $\fS_{(n)}=\fS_n$, $X_{(n)}=\{1\}$ and 
$d_{(n)}=1$. Consequently, we have
\[ \fR((n))=\{w_{(n)}\} \qquad  \mbox{and}\qquad \fN^{(n)}= 
\langle C_{w_{(n)}}\rangle_A.\]
On the other hand, we also have $Z_{w_{(n)}}=C_{w_{(n)}}$. Hence the 
desired statement holds in this case. Now consider an arbitrary 
partition $\lambda \neq (n)$ and assume that the desired statements hold
for all elements in sets $\fR(\nu)$ where $\lambda\triangleleft \nu$. Let 
$1\leq i,j \leq d_\lambda$ be such that $w=\bw_\lambda(i,j)$. By
Lemma~\ref{mainlem1}, we can write 
\[Z_w=\sum_{z\in\fS_n} f_{z}\,C_z\qquad\mbox{where $f_{z}\in A$}.\]

\medskip
\noindent {\em Claim (1).} {\em We have $f_z \in \Z$ for all $z\in \fS_n$ 
such that $\lambda \not{\!\triangleleft} \;\lambda_z$}.
\medskip

This is seen as follows. Let $\cA$ be the collection of all $z\in \fS_n$ 
such that $f_z\neq 0$ and $\lambda \not{\!\triangleleft} \;\lambda_z$. 
If $\cA=\varnothing$, there is nothing to be proved. Now assume
that $\cA\neq \varnothing$. Since $Z_w=\overline{Z}_w$, we have $f_z=
\overline{f}_z$ for all $z\in \fS_n$. Hence it will be sufficient to show
that the non-negative number 
\[m:=\min \{i \geq 0 \mid v^if_z \in {\Z}[v] \mbox{ for all $z\in \cA$}\}\]
is actually equal to $0$.
Let $z'\in \cA$ be such that $v^mf_{z'}\in {\Z}[v]$ has a non-zero constant 
term. Let $\mu\in\Lambda_n$ be such that $z'\in\fR(\mu)$ and write $z'=
\bw_\mu(k,l)$ where $1\leq k,l\leq d_\mu$. Note that, by the definition of 
$\cA$, we have $\lambda \not{\!\triangleleft}\;\mu$. Now we obtain the
following identity:
\[ v^{m+a(\mu)}\, \fX_\mu^{kl}(Z_w)=\sum_{z\in\fS_n} 
v^m\,f_z\,\Bigl(v^{a(\mu)}\,\fX_\mu^{kl}(C_z)\Bigr).\]
Let $z\in \fS_n$ be such that $z\in \fR(\nu)$ where $\lambda\triangleleft\nu$.
By induction, $C_z\in \fN^\nu$. Now Corollary~\ref{nonzero} shows that, 
if $\fX_\mu(C_z)\neq 0$, then $\nu \trianglelefteq \mu$ and so $\lambda 
\triangleleft \mu$, a contradiction.  Hence, the above sum need only be 
extended over all elements $z\in \cA$. But then we have $v^mf_z\in {\Z}[v]$ 
and so 
\[ v^{m+a(\mu)}\,\fX_\mu^{kl}(Z_w) \equiv \sum_{z\in\cA} v^m\,f_z\,
\varepsilon_z\,c_{z,\mu}^{kl} \quad \bmod \fp.\]
Now, we have $c_{z,\mu}^{kl}=0$ unless $z=z'$; see Lemma~\ref{lem01}. 
Hence we obtain 
\[ v^{m+a(\mu)}\, \fX_\mu^{kl}(Z_w) \equiv v^m\,f_{z'}\, \varepsilon_{z'}\,
c_{z',\mu}^{kl} \equiv \pm v^m\, f_{z'} \quad \bmod \fp.\]
Since $v^mf_{z'}$ has a non-zero constant term, we conclude that the above
expression is not congruent to $0$ modulo $\fp$. In particular, the left
hand side is non-zero. Since $Z_w\in \fN^\lambda$, we can now deduce
that $\lambda \trianglelefteq \mu$; see Corollary~\ref{nonzero}. Since we 
also have $\lambda \not{\!\triangleleft} \,\mu$, we conclude that 
$\lambda=\mu$. Thus, we have reached the conclusion  that 
\[ v^{m}\, \Bigl(v^{a(\lambda)}\, \fX_\lambda^{kl}(Z_w)\Bigr) \not\equiv 
0 \quad \bmod \fp.\]
Using Lemma~\ref{mainlem2}, we now see that $m=0$. Thus, (1) is proved.

\medskip
\noindent {\em Claim (2).} {\em We have $f_z=0$ unless $\lambda 
\trianglelefteq \lambda_z$}.
\medskip

To see this, let $\mu \in \Lambda_n$ be such that $\lambda 
\not{\!\!\trianglelefteq} \,\mu$. Assume, if possible, that there exists 
some $z'\in \fR(\mu)$ such that $f_{z'}\neq 0$. Let $1\leq k,l\leq d_\mu$ 
be such that $z'=\bw_\mu(k,l) \in \fR(\mu)$. Now, as above, we see that 
\[ v^{a(\mu)}\, \fX_\mu^{kl}(Z_w)=\sum_{z} f_{z}\,\Bigl(v^{a(\mu)}\,
\fX_\mu^{kl}(C_{z})\Bigr)\]
where the sum need only be extended over all $z\in \fS_n$ such that
$\lambda\not{\!\!\trianglelefteq} \,\lambda_z$. By (1), we have $f_{z}
\in \Z$ for all such elements~$z$. Hence we obtain
\[ v^{a(\mu)}\, \fX_\mu^{kl}(Z_w) \in \cO \quad \mbox{and}\quad
 v^{a(\mu)}\, \fX_\mu^{kl}(Z_w) \equiv \pm f_{z'} \not\equiv 0 
\quad \bmod \fp,\]
since $c_{z',\mu}^{kl}=\pm 1$ and $c_{z,\mu}^{kl}=0$ if $z\neq z'$.
In particular, we can now conclude that $\fX_\mu(Z_w)\neq 0$. But then 
Lemma~\ref{nonzero} shows that $\lambda \trianglelefteq \mu$, a 
contradiction. Thus, (2) is proved.

\medskip
\noindent {\em Claim (3).} {\em We have $\pm C_w\in Z_w +
\hat{\fN}^\lambda$}.
\medskip

This is seen as follows. By (1) and (2), we can write
\[ Z_w\equiv \sum_{z\in \fR(\lambda)} f_z\, C_z\quad \bmod \quad 
\langle C_y \mid y \in \fS_n,\lambda\triangleleft \lambda_y\rangle_A,\]
where $f_z \in \Z$ for all $z\in \fR(\lambda)$. By induction, we have 
$C_y \in \hat{\fN}^\lambda$ for any $y\in \fS_n$ such that $\lambda
\triangleleft \lambda_y$. Thus, we have 
\[ Z_w\equiv \sum_{z\in \fR(\lambda)} f_z\, C_z \quad \bmod  
\hat{\fN}^\lambda \qquad \mbox{where $f_z\in \Z$}.\]
Now fix $z_0\in \fR(\lambda)$ and write $z_0=\bw_\lambda(k,l)$ where
$1\leq k,l \leq d_\lambda$. In order to determine $f_{z_0}$, we  multiply 
the above relation by $v^{a(\lambda)}$ and apply $\fX_\lambda^{kl}$. By 
Corollary~\ref{nonzero}, we have $\fX_\lambda(h)=0$ for all $h \in 
\hat{\fN}^\lambda$. Hence we obtain 
\[ v^{a(\lambda)}\, \fX_\lambda^{kl}(Z_w)= \sum_{z\in \fR(\lambda)}
f_{z}\, \Bigl(v^{a(\lambda)}\,\fX_{\lambda}^{kl}(C_{z})\Bigr).\]
Now let us take constant terms. On the left, we obtain $\pm 1$ if $k=i$ 
and $l=j$, and $0$ otherwise; see Lemma~\ref{mainlem2}. On the right, we 
obtain 
\[ \sum_{z\in \fR(\lambda)} f_{z}\, \varepsilon_z\, c_{z,\lambda}^{kl}=
\pm f_{z_0}.\]
Hence, we obtain $f_{z_0}=\pm 1$ if $z_0=w$ and $f_{z_0}=0$, otherwise; 
thus, (3) is proved and the proof of (a) is complete. 

(b) Let $\fM^\lambda$ be the $A$-submodule of $\bH$ defined by the right 
hand side of the desired identity in (b). Let $\mu \in \Lambda_n$ be
such that $\lambda \trianglelefteq \mu$ and $w \in \fR(\mu)$. Since 
$Z_w \in \fN^\mu$, we can deduce from (a) that $C_w \in \fN^\mu \subseteq
\fN^\lambda$. Hence, we have 
\[ \fM^\lambda \subseteq \fN^\lambda.\]
Now note that $|\bT(\lambda)|=d_\lambda$; see Corollary~\ref{ident}. Hence, 
both $\fN^\lambda$ and $\fM^\lambda$ are free $A$-modules of the same
rank, namely,
\[\sum_{\atop{\mu \in \Lambda_n}{\lambda \trianglelefteq \mu}} 
|\fR(\lambda)|=\sum_{\atop{\mu \in \Lambda_n}{\lambda \trianglelefteq \mu}} 
|d_\lambda|^2=\sum_{\atop{\mu \in \Lambda_n}{\lambda \trianglelefteq \mu}} 
|\bT(\lambda)|^2.\]
Consequently, we have $K_0 \otimes_A \fN^\lambda=K_0\otimes_A \fM^\lambda$,
where $K_0$ is the field of fractions of $A$. So there exists some 
$0\neq f \in A$ such that 
\[ f\fN^\lambda \subseteq \fM^\lambda \subseteq \fN^\lambda, \]
that is, $\fN^\lambda/\fM^\lambda$ is contained in the torsion part of
$\bH/\fM^\lambda$. But, since the generators of $\fM^\lambda$ can be 
extended to an $A$-basis of $\bH$, the quotient $\bH/\fM^\lambda$ is a 
free $A$-module. Hence, we also have $\fN^\lambda\subseteq \fM^\lambda$, 
as desired.
\end{proof}

\begin{cor} \label{bridge} Let $\lambda \in \Lambda_n$ and $\fs,\ft
\in \bT(\lambda)$. Then there exists a unique element $w\in \fR(\lambda)$
such that 
\begin{align*}
\ty_{\fs\ft}=\eta_w\, C_w&+ \mbox{$v{\Z}[v]$-combination of 
elements $C_y$ where $\lambda=\lambda_y$}\\ &+\mbox{$A$-combination of 
elements $C_y$ where $\lambda\triangleleft \lambda_y$}.
\end{align*}
We have $w=\bw_\lambda(i,j)$, where $1\leq i,j\leq d_\lambda$ are such 
that $d(\fs)=x_i$ and $d(\ft)=x_j$.
\end{cor}

Applying the ring involution $j \colon \bH \rightarrow \bH$ and using
Corollary~\ref{murphib}, we also obtain a relation with the original Murphy 
basis:
\begin{align*}
v^{-l(w_\lambda)}x_{\fs\ft}=\eta_w'C_w'&+ 
\mbox{\em $v^{-1}{\Z}[v^{-1}]$-combination of elements $C_y'$ where 
$\lambda=\lambda_y$}\\ &+\mbox{\em $A$-combination of elements $C_y'$ 
where $\lambda\triangleleft \lambda_y$},
\end{align*}
where $\eta_w'=\varepsilon_{d(\fs)}\,\varepsilon_{d(\ft)}\,
\varepsilon_{w_\lambda}\, \varepsilon_w\, \eta_w=\pm 1$.

\begin{proof} We begin with the following preliminary remarks.
Let $x \in X_\lambda$. In the proof of Lemma~\ref{mainlem1}, we have
used \cite[Prop.~3.3]{myind} to write $C_{xw_\lambda}$ as a linear
combination of terms $T_yC_{w_\lambda}$ where $y \in X_\lambda$. But,
the result in [{\em loc.\ cit.}] actually yields something stronger,
namely, we have 
\[ C_{xw_\lambda}=T_{x}C_{w_\lambda}+\sum_{\atop{y\in X_\lambda}{y<x}}
a_{xy}\, T_yC_{w_\lambda} \quad \mbox{where}\quad a_{xy}\in v{\Z}[v].\]
(Note that, in \cite{myind}, we work with the $C'$-basis; applying
the ring involution $j \colon \bH \rightarrow \bH$ yields the above
reformulation.) Inverting the above relations, we also obtain that
\[ T_xC_{w_\lambda}=C_{xw_\lambda}+\sum_{\atop{y\in X_\lambda}{y<x}}
b_{xy}\, C_{yw_\lambda} \quad \mbox{where}\quad b_{xy}\in v{\Z}[v].\]
Now assume that $x=x_i$ where $1\leq i\leq d_\lambda$ is such that 
$d(\fs)=x_i$. Let $y \in X_\lambda$ be such that $y<x$ and $b_{xy}\neq 0$. 
By Theorem~\ref{ourmain}(b), we have either $yw_\lambda \in \fR(\lambda)$ or
$C_{yw_\lambda} \in \hat{\fN}^\lambda$. Hence, using Lemma~\ref{technic},
we can write 
\[ T_{x_i}C_{w_\lambda}\equiv C_{x_iw_\lambda}+\sum_{k=1}^{d_\lambda}
b_{ik}\, C_{x_kw_\lambda} \quad \bmod \hat{\fN}^\lambda\]
where $b_{ik}\in v{\Z}[v]$. Applying the anti-involution $\flat \colon
T_x \mapsto T_{x^{-1}}$, we can also write
\[ C_{w_\lambda}T_{x_j^{-1}}\equiv C_{w_\lambda x_j^{-1}}+
\sum_{l=1}^{d_\lambda} b_{jl}\, C_{w_\lambda x_l^{-1}} \quad \bmod 
\hat{\fN}^\lambda\]
where $b_{jl}\in v{\Z}[v]$ and $1\leq j \leq d_\lambda$ is such that
$d(\ft)=x_j$. Consequently, we obtain that 
\[ T_{x_i}C_{w_\lambda}^2 T_{x_j^{-1}} \equiv C_{x_iw_\lambda}\,
C_{w_\lambda x_j^{-1}}+\sum_{k,l=1}^{d_\lambda} c_{kl}\, C_{x_kw_\lambda}\, 
C_{w_\lambda x_l^{-1}} \quad \bmod \hat{\fN}^\lambda\]
where $c_{kl} \in v{\Z}[v]$. Using the identity in Lemma~\ref{cwi} and
the definition of $Z_w$, we obtain that
\[ \ty_{\fs\ft} \equiv Z_{\bw_\lambda(i,j)} +
\sum_{k,l=1}^{d_\lambda} c_{kl}\, Z_{\bw_\lambda(k,l)} \quad \bmod 
\hat{\fN}^\lambda.\]
Finally, Theorem~\ref{ourmain}(a) yields
\[\ty_{\fs\ft} \equiv \eta_{ij}\,C_{\bw_\lambda(i,j)} +
\sum_{k,l=1}^{d_\lambda} c_{kl}\, \eta_{kl}\,
C_{\bw_\lambda(k,l)} \quad \bmod \hat{\fN}^\lambda,\]
where we set $\eta_{kl}:=\eta_y$ if $y=\bw_\lambda(k,l)$. Thus, the 
desired assertion holds where $w=\bw_\lambda(i,j)$.
\end{proof}

\begin{exmp} \label{s4} Let $n=4$ and $\lambda=(3,1)$. We have $w_{(3,1)}=
s_1s_2s_1$;  the standard $(3,1)$-tableaux are given by 
\begin{center}
\begin{picture}(320,40)
\put(00,17){$\ft^{(3,1)}=$}
\put(40, 5){\line(1,0){15}}
\put(40,20){\line(1,0){45}}
\put(40,35){\line(1,0){45}}
\put(40, 5){\line(0,1){30}}
\put(55, 5){\line(0,1){30}}
\put(70,20){\line(0,1){15}}
\put(85,20){\line(0,1){15}}
\put(45,24){$1$}
\put(60,24){$2$}
\put(75,24){$3$}
\put(45,9){$4$}
\put(110,17){$s_3.\ft^{(3,1)}=$}
\put(160,5){\line(1,0){15}}
\put(160,20){\line(1,0){45}}
\put(160,35){\line(1,0){45}}
\put(160,5){\line(0,1){30}}
\put(175,5){\line(0,1){30}}
\put(190,20){\line(0,1){15}}
\put(205,20){\line(0,1){15}}
\put(165,24){$1$}
\put(180,24){$2$}
\put(195,24){$4$}
\put(165,9){$3$}
\put(230,17){$s_2s_3.\ft^{(3,1)}=$}
\put(290,5){\line(1,0){15}}
\put(290,20){\line(1,0){45}}
\put(290,35){\line(1,0){45}}
\put(290,5){\line(0,1){30}}
\put(305,5){\line(0,1){30}}
\put(320,20){\line(0,1){15}}
\put(335,20){\line(0,1){15}}
\put(295,24){$1$}
\put(310,24){$3$}
\put(325,24){$4$}
\put(295,9){$2$}
\end{picture}
\end{center}
Thus, we have $\{d\in X_{(3,1)}\mid d.\ft^{(3,1)}\in \bT(1,3)\}=\{1, s_3,
s_2s_3\}$. Given $\fs,\ft\in \bT(3,1)$, we write $\ty_{d(\fs),d(\ft)}$ 
instead of $\ty_{\fs\ft}$. With this convention, the nine Murphy basis
elements corresponding to pairs of standard $(3,1)$-tableaux are given as
follows:
\begin{align*}
\ty_{1,1} &= \underline{C}_{s_1s_2s_1},\\ 
\ty_{1,s_3}&=vC_{s_1s_2s_1}+\underline{C}_{s_1s_2s_1s_3},\\ 
\ty_{1,s_2s_3}&=vC_{s_1s_2s_1s_3}+\underline{C}_{s_1s_2s_1s_3s_2}\\
\ty_{s_3,1}&=vC_{s_1s_2s_1}+\underline{C}_{s_1s_3s_2s_1},\\ 
\ty_{s_3,s_3}&=v^2C_{s_1s_2s_1}+vC_{s_1s_2s_1s_3}+vC_{s_1s_3s_2s_1}+
\underline{C}_{s_1s_2s_3s_2s_1},\\
\ty_{s_3,s_2s_3}&=v^2C_{s_1s_2s_1s_3}+\underline{C}_{s_1s_2s_3s_2}+
vC_{s_1s_2s_1s_3s_2}\\
&\qquad\qquad\qquad\qquad+ vC_{s_1s_2s_3s_2s_1}+C_{s_1s_2s_1s_3s_2s_1},\\
\ty_{s_2s_3,1}&=vC_{s_1s_3s_2s_1}+\underline{C}_{s_2s_1s_3s_2s_1},\\
\ty_{s_2s_3,s_3}&=v^2C_{s_1s_3s_2s_1}+\underline{C}_{s_2s_3s_2s_1}+
vC_{s_1s_2s_3s_2s_1}\\ &\qquad\qquad\qquad\qquad
+ vC_{s_2s_1s_3s_2s_1}+ C_{s_1s_2s_1s_3s_2s_1},\\
\ty_{s_2s_3,s_2s_3}&=\underline{C}_{s_2s_3s_2}+vC_{s_1s_2s_3s_2}+
vC_{s_2s_3s_2s_1}\\ &\qquad\qquad\qquad\qquad +v^2C_{s_1s_2s_3s_2s_1}+
(v-v^{-1}) C_{s_1s_2s_1s_3s_2s_1}
\end{align*}
In each case, we have underlined the Kazhdan--Lusztig basis element 
which corresponds to the Murphy basis element as in Corollary~\ref{bridge}.
These examples show, in particular, that the coefficient of $C_y$,
where $\lambda \triangleleft \lambda_y$, may involve negative powers
of~$v$.
\end{exmp}

\section{Applications to the Kazhdan--Lusztig cells in $\fS_n$} 
\label{sec-app}

We keep the setting of the previous sections. In particular, recall the 
partition 
\[ \fS_n =\coprod_{\lambda\in \Lambda_n} \fR(\lambda) \qquad \mbox{where}
\qquad w_\lambda \in \fR(\lambda),\]
see Example~\ref{symgroup}. By Remark~\ref{incells}, we already know
that each set $\fR(\lambda)$ is contained in a two-sided cell of $\fS_n$.
Now we can prove the following result. 

\begin{thm} \label{ourmain1} Let  $\lambda,\mu \in \Lambda_n$ and
$x \in \fR(\lambda)$, $y \in \fR(\mu)$. Then we have 
\[ x \leq_{\cLR} y\qquad\Leftrightarrow\qquad\mu\trianglelefteq \lambda.\]
In particular, $w_\lambda \leq_{\cLR} w_\mu$ if and only if $\mu
\trianglelefteq \lambda$. The sets $\fR(\lambda)$, $\lambda 
\in \Lambda_n$, are precisely the two-sided cells of $\fS_n$.
\end{thm}

\begin{proof} If $\mu\trianglelefteq \lambda$, then $w_\lambda \leq_{\cLR}
w_\mu$ by Remark~\ref{kostka}. Since each of the sets $\fR(\lambda)$ and
$\fR(\mu)$ is contained in a two-sided cell, we conclude that $x \leq_{\cLR}
y$. Conversely, assume that $x \leq_{\cLR} y$. We may assume without
loss of generality that $x \leftarrow_{\cL} y$ or $x \leftarrow_{\cR} y$.
Let us first assume that $x\leftarrow_{\cL} y$. This means that exists 
some $h \in \bH$ such that $C_x$ occurs in $hC_y$ (expressed in the
$C$-basis). Since $y \in \cR(\mu)$, we have $C_y \in \fN^\mu$ by
Theorem~\ref{ourmain}(a). Hence we also have $hC_y \in \fN^\mu$ since
$\fN^\mu$ is a two-sided ideal. But then Theorem~\ref{ourmain}(b) shows 
that $hC_y$ can be written as an $A$-linear combination of elements
$C_z$ where $\mu\trianglelefteq \lambda_z$. Since $C_x$ occurs in that
linear combination, we conclude that $\mu \trianglelefteq \lambda$, as 
desired. The argument in the case where $x\leftarrow_{\cR} y$ is 
completely analogous. The statement concerning the two-sided cells is 
now clear.
\end{proof} 

The equivalence ``$w_\lambda \leq_{\cLR} w_\mu \Leftrightarrow \mu 
\trianglelefteq \lambda$'' is not a new result. However, as far as we are
aware of, all the previously known proofs for the implication 
``$\Rightarrow$'' rely on the geometric interpretation of the 
Kazhdan--Lusztig basis and the resulting ``positivity properties''. 
See, for example, Du--Parshall--Scott \cite[2.13.1]{DPS}, 
where this is deduced from a result of Lusztig--Xi \cite[3.2]{LX}.

\begin{lem} \label{lem51} Let $\lambda \in \fR(\lambda)$ and $1\leq
j \leq d_\lambda$. Then 
\[ V_j^\lambda:=\langle C_{\bw_\lambda(i,j)} +\hat{\fN}^\lambda \mid 1\leq 
i \leq d_\lambda\rangle_A \subseteq \fN^\lambda/\hat{\fN}^\lambda\]
is a left $\bH$-module, free over $A$ of rank $d_\lambda$. We have
\[ \fN^\lambda/\hat{\fN}^\lambda=V_1^\lambda \oplus \cdots \oplus 
V_{d_\lambda}^\lambda.\]
\end{lem}

\begin{proof} First note that, by Theorem~\ref{ourmain}, $\fN^\lambda/
\hat{\fN}^\lambda$ is free over $A$, with a basis given by the residue
classes of the elements $C_w$ where $w \in \fR(\lambda)$. This already
shows that each $V_j^\lambda$ is free over $A$ of rank $d_\lambda$, and that
$\fN^\lambda/\hat{\fN}^\lambda$ is the direct sum of all $V_j^\lambda$.
It remains to prove that $V_j^\lambda$ is a left $\bH$-module. For
this purpose, we must prove that, for any $h \in \bH$ and $1\leq i,k\leq 
d_\lambda$, there exist $r_h^j(i,k) \in A$ such that 
\[ hC_{\bw_\lambda(i,j)}\equiv \sum_{k=1}^{d_\lambda} r_h^j(i,k)\, 
C_{\bw_\lambda(k,j)} \quad \bmod \hat{\fN}^\lambda.\]
First we prove this in the case where $j=1$. Then 
\[ \bw_\lambda(i,1)=x_iw_\lambda \qquad \mbox{where $x_i \in X_\lambda$};\]
see Remark~\ref{choice}. Now, by Lemma~\ref{cwi}(c), the $A$-module
\[ \langle C_{xw_\lambda} \mid x\in X_\lambda\rangle_A\]
is a left ideal in $\bH$. Hence, for any $h \in \bH$, we can write
\[ hC_{x_iw_\lambda}=\sum_{x\in X_\lambda} r_h(i,x)\, C_{xw_\lambda}
\qquad \mbox{where $r_h(i,x)\in A$}.\]
Now let $x\in X_\lambda$ and $\mu \in \Lambda_n$ be such that
$xw_\lambda \in \fR(\mu)$. Then, by Lemma~\ref{technic}, we have  
either $\lambda \triangleleft \mu$ or $x=x_k$ for some $k\in
\{1,\ldots,d_\lambda\}$. By Theorem~\ref{ourmain}, the former condition 
implies that $C_{xw_\lambda}\in\hat{\fN}^\lambda$.  Consequently, we have 
\[ hC_{x_iw_\lambda}\equiv \sum_{k=1}^{d_\lambda} r_h(i,x_k)\,
C_{x_kw_\lambda} \quad \bmod \hat{\fN}^\lambda.\]
Thus, the desired assertion holds where we set $r_h^1(i,k):=r_h(i,x_k)$.

Now let $j>1$. By Theorem~\ref{ourmain}, there exist signs $\eta_{ij}
=\pm 1$ such that 
\[ C_{x_iw_\lambda}C_{w_\lambda x_j^{-1}} \equiv \eta_{ij}\, 
C_{\bw_\lambda(i,j)} \quad \bmod \hat{\fN}^\lambda \qquad \mbox{for
$1\leq i \leq d_\lambda$}.\]
We shall now set 
\[ r_h^j(i,k):=\eta_{ij} \,\eta_{kj}\, r_h^1(i,k) 
\qquad \mbox{for $1\leq i,k\leq d_\lambda$}.\]
Then, since $\hat{\fN}^\lambda$ is a two-sided ideal, we obtain 
\begin{align*}
hC_{\bw_\lambda(i,j)} &\equiv \eta_{ij}(h C_{x_iw_\lambda})
C_{w_\lambda x_j^{-1}} \\
&\equiv \eta_{ij} \sum_{k=1}^{d_\lambda} r_h^1(i,k)\, 
C_{x_kw_\lambda} C_{w_\lambda x_j^{-1}}\\
&\equiv \eta_{ij} \sum_{k=1}^{d_\lambda} \eta_{ki}\, r_h^1(i,k)\, 
C_{\bw_\lambda(k,j)}\\
&\equiv \sum_{k=1}^{d_\lambda} r_h^j(i,k)\, C_{\bw_\lambda(k,j)}
\end{align*}
as desired, where all of the above congruences are taken modulo 
$\hat{\fN}^\lambda$.
\end{proof}

\begin{thm} \label{ourmain2} Let $x,y \in \fS_n$. Then we have the
following implication:
\[x\leq_{\cL} y \quad \mbox{and} \quad x\sim_{\cLR} y \quad\Rightarrow 
\quad x \sim_{\cL} y.\]
\end{thm}

\begin{proof} Let $x,y \in \fS_n$ be such that $x\leq_{\cL} y$ and
$x \sim_{\cLR} y$. Since the relation $\leq_{\cL}$ is defined as the
transitive closure of $\leftarrow_{\cL}$, it will be sufficient to 
consider the special case where $x\leftarrow_{\cL} y$. Now, since 
$x\sim_{\cLR} y$, we have $x,y \in \fR(\lambda)$ for some $\lambda 
\in \fR(\lambda)$; see Theorem~\ref{ourmain1}. Thus, we have 
$y=\bw_\lambda(i,j)$ where $1\leq i,j\leq d_\lambda$. By Lemma~\ref{lem51},
we have 
\[ hC_y \equiv hC_{\bw_\lambda(i,j)}\equiv \sum_{k=1}^{d_\lambda} 
r_h^j(i,k)\, C_{\bw_\lambda(k,j)} \quad \bmod \hat{\fN}^\lambda,\]
where $r_h^j(i,k)\in A$. By Theorem~\ref{ourmain} and Theorem~\ref{ourmain1}, 
$\hat{\fN}^\lambda$ is spanned by elements $C_w$ where $w\leq_{\cLR} y$ 
but $w\not\sim_{\cLR} y$. Hence, our hypothesis that $C_x$ occurs in
the decomposition of $hC_y$ implies that either $x=\bw_\lambda(k,j)$
for some $1\leq k \leq d_\lambda$ or $x \not\sim_{\cLR} y$. Hence, we
must have $x=\bw_\lambda(k,j)$ and so $x\sim_{\cL} y$ by 
Remark~\ref{incells}.
\end{proof} 

Again, the above result is not new (and the conclusion is known to hold 
for more general types of Coxeter groups). But, even for the symmetric 
group $\fS_n$, it was first proved by Lusztig \cite[Lemma~4.1]{Lu0} using 
the (deep) connection between cells and primitive ideals in universal 
envelopping algebras via the main conjecture of Kazhdan--Lusztig 
\cite{KaLu}. Another proof, applicable to finite and affine Weyl groups, 
was given by Lusztig \cite{Lu1}, using the geometric interpretation of the 
Kazhdan--Lusztig basis and the resulting ``positivity properties''.

The following result was first obtained by Kazhdan--Lusztig 
\cite[Theorem~1.4]{KaLu}, as a consequence of the combinatorial 
description of the left cells in terms of the Robinson--Schensted
correspondence. (We will come back to the latter point in
Theorem~\ref{finalrs}.)

\begin{thm} \label{ourmain3} For each left cell $\fC$ of $\fS_n$, we have 
$\chi_{\fC} \in \Irr(\bH_K)$. Furthermore, two left cells $\fC,\fC'$ 
afford the same character if and only if $\fC,\fC'\subseteq \fR(\lambda)$ 
for some $\lambda \in \Lambda_n$. The total number of left cells equals
$\sum_{\lambda\in \Lambda_n} d_\lambda$.
\end{thm}

\begin{proof} For $\lambda \in \Lambda_n$ and  $1\leq j \leq d_\lambda$,
we set 
\[ {^j\cC}_\lambda:=\{\bw_\lambda(i,j) \mid 1\leq i \leq d_\lambda \}
\subseteq \fR(\lambda).\]
By Remark~\ref{incells}, the above set is contained in a left cell of
$\fS_n$. We claim that ${^j}\cC_\lambda$ is a left cell. Since $\leq_{\cL}$ 
is defined as the transitive closure of the relation $\leftarrow_{\cL}$, 
it is enough to prove the following statement. 

\medskip
{\em Let $y \in {^j\cC}_\lambda$ and $h \in \bH$. Then $hC_y$ is a linear 
combination of basis elements $C_x$ where $x\in \fC_j^\lambda$ or 
$x<_{\cLR} y$.}

\medskip
This is proved as follows. By Lemma~\ref{lem51}, we can write $hC_y$ as a
linear combination of basis elements $C_x$ where $x\in {^j}\cC_\lambda$, 
and some element of $\hat{\fN}^\lambda$. But, by Theorem~\ref{ourmain} and 
Theorem~\ref{ourmain1}, $\hat{\fN}^\lambda$ is spanned by basis 
elements $C_w$ where $w<_{\cLR} y$. This yields the above statement.
Thus, ${^j\cC}_\lambda$ is a left cell, as claimed.  Now it is clear that
we have a partition
\[ \fS_n=\coprod_{\lambda \in \Lambda_n} \coprod_{j=1}^{d_\lambda}
{^j\cC}_\lambda.\]
Hence, the sets ${^j\cC}_\lambda$ are precisely the left cells of $\fS_n$.
The remaining statements concerning the characters of the left cell 
representations now follow by a standard counting argument. Indeed, since
$\bH_K$ is semisimple, the above partition yields an isomorphism of left
$\bH_K$-modules
\[ \bH_K \cong \bigoplus_{\lambda \in \Lambda_n} \bigoplus_{j=1}^{d_\lambda}
\;[{^j\cC}_\lambda]_K.\]
Thus, $\bH_K$ (regarded as a left $\bH_K$-module) has a decomposition into
a direct sum with $\sum_\lambda d_\lambda$ terms. Now recall that the latter
number is the sum of the dimensions of all irreducible representations
of $\bH_K$ (up to equivalence). Hence, by Wedderburn's Theorem, 
$\sum_\lambda d_\lambda$ is the maximum number of terms in a direct sum 
decomposition of $\bH_K$. Consequently, each direct summand 
$[{^j\cC}_\lambda]_K$ must be a simple $\bH_K$-module. Thus, we have shown
that $[\fC]_K$ is a simple $\bH_K$-module for any left cell $\fC$. 
Finally, we claim that 
\[ \chi_{\lambda}=\mbox{character afforded by $[{^j\cC}_\lambda]_K$}.\]
Indeed, since ${^j\cC}_\lambda\subseteq\fR(\lambda)$, Proposition~\ref{mylem}
shows that $\chi_\lambda$ occurs in the character afforded by 
$[{^j\cC}_\lambda]_K$. Since the latter is irreducible, we have equality,
as desired. We conclude that two left cells afford the same character if
and only if these two left cells are both contained in  $\fR(\lambda)$ for 
some $\lambda$.
\end{proof}

Let us recall some basic facts concerning the Robinson--Schensted
correspondence. We use Knuth \cite[\S 5.1.4]{Knuth} as a reference. Recall
that, for any $\lambda \in \Lambda_n$, we denote by $\bT(\lambda)$ the set 
of all standard $\lambda$-tableaux. Then the Robinson--Schensted 
correspondence defines a bijection
\[ \coprod_{\lambda \in \Lambda_n} \bigl(\bT(\lambda) \times
\bT(\lambda)\bigr) \rightarrow \fS_n.\]
If $\fs,\ft \in \bT(\lambda)$, we denote by $\pi_\lambda(\fs,\ft)$ the
corresponding element of $\fS_n$. Given $w\in \fS_n$, the pair $(\fs,\ft)$ 
such that $w=\pi_\lambda(\fs,\ft)$ can be constructed explicitly by the 
``insertion algorithm'' \cite[p.~49]{Knuth}. The tableau $\fs$ is obtained 
by ``inserting'' the numbers from the sequence $(w.1,w.2,\ldots,w.n)$ 
into an initially empty tableau; the tableau $\ft$ ``keeps the record'' 
of the order in which the positions in $\fs$ have been filled. For example,
applying the insertion algorithm to the element $w_\lambda$ where
$\lambda \in \Lambda_n$, we obtain
\[ w_\lambda=\pi_{\lambda^*}\bigl(\ft^{\lambda^*},\ft^{\lambda^*}\bigr)\]
where $\lambda^*$ denotes the conjugate partition. 

We shall need the following property of the Robinson--Schensted
correspondence. 

\begin{prop}[Knuth] \label{knuth} Let $\lambda \in \Lambda_n$ and 
$\fs,\fs',\ft,\ft'\in \bT(\lambda)$. We set
\[ w:=\pi_\lambda(\fs,\ft) \qquad \mbox{and}\qquad w':=
\pi_\lambda(\fs',\ft').\]
Then we have $\fs=\fs'$ if and only if $w,w'$ are linked by a finite 
sequence of ``star operations'', that is, there exists a sequence $w=y_0,
y_1, \ldots,y_k=w'$ of elements $y_i\in \fS_n$ such that $y_{i-1} \approx 
y_i$ for all $i$, where $\approx$ is defined in (\ref{star}). 
\end{prop}

For the proof, see Exc.~4 in \cite[\S 5.1.4]{Knuth}. Note that
Knuth shows that we have $\fs=\fs'$ if and only if the permutations $w,w'$
can be transformed to each other by a finite sequence of so-called
``admissible transformations''. It is readily checked that the latter 
condition is equivalent to the fact that $w,w'$ are linked by a finite
sequence of ``star operations'', using the characterisation of $l(w)$ as 
the number of pairs $(i,j)$ such that $1 \leq i<j\leq n$ and $w(j)>w(i)$. 
See also \cite[\S 3.2]{Ar} for more details.

Now we can prove the following result. The statements in (a) and (b) are due
to Kazhdan--Lusztig \cite[\S 5]{KaLu}, but the proof given there is quite
sketchy. A complete, self-contained proof, based on the methods in
\cite[\S 4]{KaLu}, is given by Ariki \cite{Ar}. The proof that we give here
is different as far as the (more difficult) implications ``$\Rightarrow$'' 
are concerned. Note also that (c) is neither proved  in \cite{KaLu}  nor
in \cite{Ar}.

\begin{cor} \label{finalrs} Let $w,w'\in \fS_n$. Let $\lambda, \mu \in 
\Lambda_n$ be such that $w=\pi_\lambda(\fs,\ft)$ where $\fs,\ft\in 
\bT(\lambda)$ and $w'=\pi_\mu(\fs',\ft')$ where $\fs',\ft'\in \bT(\mu)$. 
Then the following hold:
\begin{align*}
w\sim_{\cR} w' & \quad \Leftrightarrow \quad \lambda=\mu 
\mbox{ and } \fs=\fs', \tag{a}\\ 
w\sim_{\cL} w' & \quad \Leftrightarrow \quad \lambda=\mu \mbox{ and } 
\ft=\ft', \tag{b}\\ w\sim_{\cLR} w' & \quad \Leftrightarrow 
\quad \lambda=\mu.\tag{c}
\end{align*}
In particular, the intersection of a left cell and a right cell is either
empty or a singleton set. Furthermore, for any $\lambda\in \Lambda_n$,
we have 
\[ \fR(\lambda)=\{\pi_{\lambda^*}(\fu,\fv)\mid\fu,\fv\in\bT(\lambda^*)
\},\]
where $\lambda^*$ denotes the conjugate partition.
\end{cor}

\begin{proof} First note that $d_\lambda=|\bT(\lambda)|$ for any
$\lambda \in \Lambda_n$; see Corollary~\ref{ident}. 

(a) Let $\lambda \in \Lambda_n$ and $\fs\in \bT(\lambda)$. By
Proposition~\ref{knuth} and (\ref{star})(a), the set 
\[ \bT(\lambda,\fs):=\{ \pi_\lambda(\fs,\ft) \mid \ft\in \bT(\lambda)\}\]
is contained in a right cell. Now note that the family of sets 
\[ \{\bT(\lambda,\fs) \mid \lambda \in \Lambda_n,\fs\in \bT(\lambda)\}\]
forms a partition of $\fS_n$, and that there are $\sum_\lambda d_\lambda$ 
pieces in that partition. On the other hand, by Theorem~\ref{ourmain3}, the
latter number also equals the number of all left cells of $\fS_n$, and
that number is the same as the number of all right cells. Hence each set 
$\bT(\lambda,\fs)$ must be a right cell.

(b) We have $\pi_\lambda(\fs,\ft)^{-1}=\pi_\lambda(\ft,\fs)$ for all
$\lambda \in \Lambda_n$ and $\fs,\ft\in \bT(\lambda)$; see Theorem~B in 
\cite[\S 5.1.4]{Knuth}. Hence the assertion follows from (a).

(c) Let $\lambda \in \Lambda_n$. Then, by (a) and (b), the set 
$\{\pi_\lambda(\fs,\ft)\mid \fs,\ft\in\bT(\lambda)\}$ is contained in a 
two-sided cell of $\fS_n$. Using Theorem~\ref{ourmain1}, we can now argue 
as in the proof of (a) to conclude that the above set is a two-sided cell. 

Finally, consider the statement concerning $\fR(\lambda)$. By 
Remark~\ref{choice}, we have $w_\lambda \in \fR(\lambda)$. On the
other hand, applying the ``insertion algorithm'' to $w_\lambda$, we
see that $w_\lambda=\pi_{\lambda^*}(\fu_0,\fv_0)$ for some standard 
$\lambda^*$-tableaux $\fu_0,\fv_0$. Now Theorem~\ref{ourmain1} and (c)
imply that $\fR(\lambda)=\{\pi_{\lambda^*}(\fu,\fv) \mid \fu,\fv\in 
\bT(\lambda^*)\}$.
\end{proof}

\begin{rem} \label{palli} Once the above result is established, we can also
identify the elements $x_i$ in Remark~\ref{choice}. Indeed, for $\lambda
\in \Lambda_n$, the set 
\[ {^1}\cC_\lambda=\{x_iw_\lambda \mid 1\leq i\leq d_\lambda\}\]
is the left cell containing $w_\lambda$. Then the Robinson--Schensted
correspondence shows that
\[ {^1}\cC_\lambda=\{d(\fs)w_\lambda \mid\fs \in \bT(\lambda)\};\]
see \cite[Lemma~3.3]{McPa} and the references there. We also remark here
that, in \cite[Theorem~3.5]{McPa}, McDonough and Pallikaros construct an
explicit isomorphism of $\bH$-modules from $[{^1}\cC_\lambda]_A$ onto
the Specht module $\tilde{\cS}^{\lambda^*}$, where $\lambda^*$ is
the conjugate partition.
\end{rem}

\begin{cor} \label{maincor4} Let $\lambda \in \Lambda_n$ and $\fC,\fC_1$ 
be two left cells contained in $\fR(\lambda)$. Then there is a unique 
bijection $\fC\stackrel{\sim}{\rightarrow} {\fC}_1$, $x \mapsto x_1$, 
such that the following conditions are satisfied:
\begin{align*}
x\sim_{\cR} x_1 \qquad  & \qquad \mbox{for all $x\in \fC$},\tag{a}\\
h_{w,x,y}=h_{w,x_1,y_1} &\qquad\mbox{for all $x,y\in \fC$ and 
$w \in \fS_n$}.  \tag{b}
\end{align*}
Thus, the $\bH$-modules $[\fC]_A$ and $[\fC_1]_A$ are not only isomorphic, 
but they even afford exactly the same matrix representations with respect
to the standard bases of $[\fC]_A$ and $[\fC]_A$, respectively.
\end{cor}

\begin{proof} This follows by the argument in \cite[\S 5]{KaLu}. Indeed,
by Corollary~\ref{finalrs} and Proposition~\ref{knuth}, $\fC$ and $\fC_1$ 
can be linked by a finite sequence of star operations. Hence the 
assertions follow from (\ref{star}).
\end{proof}

Finally, we turn to the properties (P1)--(P15) in (\ref{conj}). 
Recall that, for any $w \in \fS_n$, we defined $\alpha_w:=a(\lambda)=
l(w_\lambda)$, where $\lambda \in \Lambda_n$ is such that $w \in
\fR(\lambda)$; see Definition~\ref{maindef}.

\begin{lem} \label{maincor2} Let $x,y\in \fS_n$ and assume that $x\leq_{\cLR}
y$. Then we have $\alpha_y\leq \alpha_x$, with equality only if 
$x\sim_{\cLR} y$. 
\end{lem}

\begin{proof} Let $\lambda,\mu\in\Lambda_n$ be such that $x\in \fR(\lambda)$,
$y \in \fR(\mu)$. By Corollary~\ref{ourmain1}, we have $\mu
\trianglelefteq \lambda$. This implies $a(\mu)\leq a(\lambda)$ with equality 
only if $\lambda=\mu$; see, for example, \cite[(5.4.2) and 
Exc.~5.6]{ourbuch}. Thus, we have $\alpha_y\leq \alpha_x$, with 
equality only if $\lambda=\mu$. Finally, by Theorem~\ref{ourmain1}, we 
have $x\sim_{\cLR} y$ if $\lambda=\mu$. 
\end{proof} 

\begin{thm} \label{p114} We have $\ba(w)=\alpha_w$ for any $w\in \fS_n$.
Furthermore, the properties (P1)--(P15) in (\ref{conj}) hold, where the 
set of ``distinguished involutions'' is given by 
\[ \cD=\{z\in \fS_n\mid z^2=1\}.\]
The constants $\gamma_{x,y,z}$ are given as follows. Let $x,y,z\in \fS_n$
and $d\in \cD$ be such that $y^{-1}\sim_{\cL} d$. Then we have 
\[ \gamma_{x,y,z}=\left\{\begin{array}{cl} n_d & \qquad \mbox{if
$x\sim_{\cL} y^{-1}$, $y\sim_{\cL} z^{-1}$, $z\sim_{\cL} x^{-1}$},\\
0 & \qquad \mbox{otherwise}.\end{array}\right.\] 
Furthermore, $n_d=1$ for any $d\in \cD$. Lusztig's ring $J$, see
(\ref{conj}), is isomorphic to the direct sum of the matrix rings 
$M_{d_\lambda}(\Z)$ ($\lambda \in \Lambda_n$).
\end{thm}

\begin{proof} First we determine $\ba(w)$ and show that (P1)--(P15) hold. 
Since condition ($*$) in Lemma~\ref{lem01} is satisfied, we can apply the 
results in \cite[\S 4]{GeIa05}. By Lemma~\ref{maincor2}, the hypotheses 
of \cite[Lemma~4.4]{GeIa05} are satisfied. Hence we have $\ba(w)=\alpha_w$ 
for any $w\in \fS_n$. Once the identity $\ba(w)=\alpha_w$ is established, 
Lemma~\ref{maincor2} shows that (P4) and (P11) hold. 

But then we can also apply \cite[Lemma~4.6]{GeIa05} and this yields (P1), 
(P2), (P3), (P5), (P6), (P7), (P8) and  (P14). Since all of the above 
properties hold for any parabolic subgroup of $\fS_n$, (P12) also holds;
see the argument in \cite[14.12]{Lusztig03}. By Theorem~\ref{ourmain3}, 
we have $\chi_{\fC}\in\Irr(\bH_K)$ for any left cell $\fC$. Hence 
\cite[Lemma~4.8]{GeIa05} yields that (P13) holds and that $\cD=\{z\in 
\fS_n\mid z^2=1\}$. By \cite[Remark~4.10]{GeIa05}, we have 
\[ \gamma_{x,y,z}=\left\{\begin{array}{cl} \pm 1 & \qquad \mbox{if
$x\sim_{\cL} y^{-1}$, $y\sim_{\cL} z^{-1}$, $z\sim_{\cL} x^{-1}$},\\
0 & \qquad \mbox{otherwise}.\end{array}\right.\] 
The statement about the structure of $J$ is contained in 
\cite[Lemma~4.9]{GeIa05}.

Next we show that (P9), (P10) hold. Now, by \cite[14.10]{Lusztig03}, 
property (P10) is a formal consequence of (P9). To prove (P9), let 
$x,y\in \fS_n$ be such that $x\leq_{\cL} y $ and $\ba(x)=\ba(y)$. By 
(P11), we obtain $x\sim_{\cLR} y$. Hence Theorem~\ref{ourmain2}
implies $x\sim_{\cL} y$, as desired. 

It remains to consider (P15). Since (P4), (P9), (P10) are already known
to hold, (P15) can be reformulated as explained in \cite[14.15]{Lusztig03}.
But then we can argue as in \cite[15.7]{Lusztig03} to conclude that (P15)
hold. (Note that we are in the case of equal parameters; the ``positivity
properties'' are not required in \cite[15.7]{Lusztig03}.) Thus, we have 
proved (P1)--(P15) for $W=\fS_n$. 

Next, let $x,y,z\in \fS_n$ and $d\in \cD$ be such that $x\sim_{\cL} y^{-1}$, 
$y\sim_{\cL} z^{-1}$, $z\sim_{\cL} x^{-1}$ and $y^{-1}\sim_{\cL} d$. We 
show that $\gamma_{x,y,z}=n_d$. Now, by Theorem~\ref{ourmain1}, we have 
$x,y,z\in \fR(\lambda)$ for some $\lambda \in \fR(\lambda)$. Furthermore, 
there exist $1\leq i,j \leq d_\lambda$ such that 
\[x=\bw_{\lambda}(i,j),\quad y=\bw_\lambda(j,k),\quad z=\bw_\lambda(k,i),
\quad d=\bw_\lambda(j,j).\]
(See Remark~\ref{incells} and note that, by the proof of 
Theorem~\ref{ourmain3}, the sets ${^j}\cC_\lambda$ defined there are 
precisely the left cells contained in $\fR(\lambda)$.) By 
Remark~\ref{incells}, we also have $x\sim_{\cL} d$. So (P5), (P7), (P13) 
show that
\[\gamma_{x,d,x^{-1}}=\gamma_{x^{-1},x,d}=n_d=\pm 1.\]
Now, since $d=d^{-1}\sim_{\cR} y$, $x\sim_{\cR} z^{-1}$ and $y\sim_{\cL}
z^{-1}$, we have 
\[ h_{x,d,x}=h_{x,y,z^{-1}}; \qquad \mbox{see Corollary~\ref{maincor4}}.\]
By (P4), we have $\ba(x)=\ba(z^{-1})$. Hence the above identity implies that
\begin{align*}
\gamma_{x,y,z} &= \mbox{ constant term of $v^{\ba(z^{-1})} h_{x,y,z^{-1}}$}
\\ &= \mbox{ constant term of $v^{\ba(x)}h_{x,d,x}$}\\
&=\gamma_{x,d,x^{-1}}=n_d,
\end{align*}
as desired. 

Finally, we show that $n_d=1$ for any $d\in \cD$. Let $\lambda \in
\Lambda_n$ be such that $d\in \fR(\lambda)$ and write $d=\bw_\lambda(j,j)$
where $1\leq j \leq d_\lambda$. We consider the element $x_jw_\lambda=
\bw_\lambda(j,1)$ where $x_j \in X_\lambda$; see Remark~\ref{choice}. 
Using the description of $C_{w_\lambda}$ in Lemma~\ref{cwi} and the 
multiplication rules in (\ref{multrule}), one easily shows that 
\[C_{x_j w_\lambda}C_{w_\lambda} =\varepsilon_{w_\lambda} 
v^{-l(w_\lambda)} \,P_\lambda \,C_{x_jw_\lambda}.\]
Thus, we have 
\[h_{x_jw_\lambda,w_\lambda,x_jw_\lambda}=v^{-l(w_\lambda)}\,P_\lambda.\]
Now, by Remark~\ref{incells}, we have $x_jw_\lambda \sim_{\cR} d$. 
Furthermore, since $w_\lambda x_j^{-1}= (x_jw_\lambda)^{-1}=
\bw_\lambda(1,j)$, we have $w_\lambda \sim_{\cR} w_\lambda x_j^{-1}
\sim_{\cL} d$. Hence Corollary~\ref{maincor4} implies that 
\[ h_{x_jw_\lambda,w_\lambda x_j^{-1},d}=h_{x_jw_\lambda,w_\lambda,
x_jw_\lambda }=v^{-l(w_\lambda)}\,P_\lambda.\]
By (P4), we have $\ba(d)=\ba(x_jw_\lambda)$ and this equals $a(\lambda)=
l(w_\lambda)$, as we have seen at the beginning of the proof. Since 
$P_\lambda$ is a polynomial with constant term $1$, we now conclude that 
\[ \gamma_{x_jw_\lambda,w_\lambda x_j^{-1},d}=\gamma_{x_jw_\lambda,
w_\lambda, w_\lambda x_j^{-1}}=1.\]
By (P5), the left hand side equals $n_d$. Thus, $n_d=1$. 
\end{proof}

We can now also determine the sign in Theorem~\ref{ourmain} and 
Corollary~\ref{bridge}. 

\begin{cor} \label{sign} We have $\eta_w=1$ for all $w\in \fS_n$.
\end{cor}

\begin{proof} Let $\lambda \in \Lambda_n$ be such that $w\in \fR(\lambda)$
and write $w=\bw_\lambda(i,j)$ where $1\leq i,j\leq d_\lambda$. Then we have 
\[ \varepsilon_{w_\lambda}\,P_\lambda Z_w=v^{l(w_\lambda)}\, 
C_{x_iw_\lambda}\, C_{w_\lambda x_j^{-1}}=\sum_{z\in \fS_n}
\varepsilon_{x_i}\, \varepsilon_{x_j} \,\varepsilon_z\, v^{l(w_\lambda)}\,
h_{x_iw_\lambda, w_\lambda x_j^{-1},z}\, C_z.\] 
On the other hand, by Theorem~\ref{ourmain}(a), we have $Z_w \equiv 
\eta_w\,C_w \bmod \hat{\fN}^\lambda$ where $\eta_w=\pm 1$. Using 
Theorem~\ref{ourmain}(b), we conclude that
\[ P_\lambda\, \varepsilon_{w_\lambda}\,\eta_w=\varepsilon_{x_i}\, 
\varepsilon_{x_j}\, \varepsilon_w\,v^{l(w_\lambda)}\,
h_{x_iw_\lambda,w_\lambda x_j^{-1},w}.\]
By Theorem~\ref{p114}, we have $\ba(w)=a(\lambda)=l(w_\lambda)$. Hence,
taking constant terms on both sides of the  above identity, we obtain
\[ \eta_w= \varepsilon_{x_i}\, \varepsilon_{x_j}\,\varepsilon_w\,
\varepsilon_{w_\lambda}\,\gamma_{x_iw_\lambda,w_\lambda x_j^{-1},w^{-1}}.\]
Now, in Theorem~\ref{p114}, we also have seen that the constants
$\gamma_{x,y,z}$ are either $0$ or $1$. This yields that
\[ \eta_w=\varepsilon_{x_i}\, \varepsilon_{x_j}\,\varepsilon_{w_\lambda}\,
\varepsilon_w.\]
Now let $\fC$ be the left cell containing $w_\lambda$ and $\fC_1$ be
the left cell containing $w_\lambda x_j^{-1}$. Then we have a bijection
$\fC \stackrel{\sim}{\rightarrow} \fC_1$ as in Corollary~\ref{maincor4}.
Under this bijection, $w_\lambda \in \fC$ corresponds to $w_\lambda x_j^{-1}
\in \fC_1$, and $x_iw_\lambda\in \fC$ corresponds to $w\in \fC_1$. Since
that bijection is the composition of a finite number of star operations,
a correspondingly repeated application of the relation in (\ref{star})(d) 
yields 
\[ l(w_\lambda)+l(x_iw_\lambda) \equiv l(w_\lambda x_j^{-1})+l(w)
\quad \bmod 2.\]
Consequently, we have $\eta_w=1$, as desired.
\end{proof}

\medskip
\noindent {\bf Acknowledgements.} This paper was written while the author 
enjoyed the hospitality of the Bernoulli center at the EPFL Lausanne 
(Switzerland), in the framework of the research program ``Group 
representation theory'' from january to june 2005.



\begin{thebibliography}{131}

\bibitem{Ar}
{\sc S.~Ariki}, {Robinson-Schensted correspondence and left cells}.
Combinatorial methods in representation theory (Kyoto, 1998), 1--20,
Adv. Stud. Pure Math., {\bf 28}, Kinokuniya, Tokyo, 2000.

\bibitem{DiJa1}
{\sc R.~Dipper and G.~D.~James}, Representations of Hecke algebras of 
general linear groups, {Proc.\ London Math.\ Soc.} {\bf 52} (1986), 20--52.

\bibitem{DiJa2}
{\sc R.~Dipper and G.~D.~James}, Blocks and idempotents of Hecke algebras 
of general linear groups, {Proc.\ London Math.\ Soc.} {\bf 54} (1989), 57--82.

\bibitem{DPS}
{\sc J.~Du, B.~Parshall, and L.~Scott}, Cells and $q$-Schur algebras,
Transformation Groups {\bf 3} (1998), 33--49.

\bibitem{my02}
{\sc M.~Geck}, Constructible characters, leading coefficients and left 
cells for finite Coxeter groups with unequal parameters, Represent. 
Theory {\bf 6} (2002), 1--30 (electronic).

\bibitem{myind}
{\sc M.~Geck}, On the induction of Kazhdan--Lusztig cells, Bull. London
Math. Soc. {\bf 35} (2003), 608--614.

\bibitem{GeIa05}
{\sc M.~Geck and L.~Iancu},  Lusztig's $a$-function in type $B_n$ in the
asymptotic case, preprint (april 2005); available at 
{\tt http://arXiv.org/math.RT/0504213}.

\bibitem{ourbuch}
{\sc M.~Geck and G.~Pfeiffer}, Characters of finite Coxeter
groups and Iwahori--Hecke algebras, London Math. Soc. Monographs,
New Series {\bf 21}, Oxford University Press, 2000.

\bibitem{Hoefs}
{\sc P.~N. Hoefsmit}, Representations of {H}ecke algebras of finite groups
with {BN} pairs of classical type, Ph.D. thesis, University of British
Columbia, Vancouver, 1974.

\bibitem{KaLu}
{\sc D.~Kazhdan and G.~Lusztig}, Representations of Coxeter groups and
Hecke algebras, Invent. Math. {\bf 53} (1979), 165--184.

\bibitem{Knuth}
{\sc D.~E. Knuth}, The art of computer programming, volume 3: Sorting and
Searching, Addison-Wesley, Second Edition, 1998.  

\bibitem{Lu0}
{\sc G.~Lusztig}, {On a theorem of {B}enson and {C}urtis},
J. Algebra {\bf 71} (1981), 490--498.
                                                                   
\bibitem{LuBook}
{\sc G.~Lusztig}, {\em Characters of reductive groups over a finite field},
Annals Math.\ Studies 107 (Princeton University Press, 1984).

\bibitem{Lu1}
{\sc G.~Lusztig}, Cells in affine Weyl groups, Advanced Studies in Pure
Math. {\bf 6}, Algebraic groups and related topics, Kinokuniya and
North--Holland, 1985, 255--287.

\bibitem{LX}
{\sc G.~Lusztig and N.~Xi}, {Canonical left cells in affine Weyl groups},
Advances in Math. {\bf 72} (1988), 284--288.

\bibitem{Lusztig90b} 
{\sc G.~Lusztig}, Intersection cohomology methods in representation theory,
 Proceedings of the {I}nternational {C}ongress of {M}athematics, {K}yoto, 
{J}apan, 1990 (I.~Satake, ed.), Springer-Verlag, 1991, pp.~155--174.

\bibitem{Lusztig03}
{\sc G.~Lusztig}, Hecke algebras with unequal parameters, CRM Monographs
Ser.~{\bf 18}, Amer. Math. Soc., Providence, RI, 2003.

\bibitem{Mathas99}
{\sc A.~Mathas}, Iwahori--{H}ecke algebras and {S}chur algebras of the
symmetric group, Univ. Lecture Ser., vol.~15, Amer. Math. Soc., 
Providence, RI, 1999.

\bibitem{McPa}
{\sc T.~P.~McDonough and C.~A.~Pallikaros}, On relations between the
classical and the Kazhdan--Lusztig representations of symmetric groups
and associated Hecke algebras, preprint (2004).

\bibitem{Mu1}
{\sc G.~E.~Murphy}, On the representation theory of the symmetric
groups and associated Hecke algebras, J. Algebra {\bf 152} (1992),
492--513.

\bibitem{Mu2}
{\sc G.~E.~Murphy}, The representations of Hecke algebras of type
$A_n$, J. Algebra {\bf 173} (1995), 97--121.

\bibitem{Spr}
{\sc T. A.~Springer}, {Quelques applications de la cohomologie d'intersection},
S\'{e}minare Bourbaki (1981/82), exp. 589, Ast\'erisque {\bf 92--93} (1982).
\end{thebibliography}
\end{document}